\documentclass{amsart}

\title{$G$-bundles on Abelian surfaces, hyperk\"ahler manifolds, and
stringy Hodge numbers.}
\author{Jim Bryan\\
Ron Donagi\\
Naichung Conan Leung}
\date{\today}
\address{
Department of Mathematics\\
Tulane University\\
6823 St. Charles Ave.\\
New Orleans, LA 70118}
\address{School of Mathematics\\
University of Minnesota\\
Minneapolis, MN 55455}
\address{Department of Mathematics\\
University of Pennsylvania\\
Philadelphia, PA 19104-6395}

\usepackage{diagrams}
\usepackage{verbatim}
\usepackage{eepic,epic}

\usepackage{amsmath,amsthm,amsfonts, amssymb}

\newtheorem{thm}{Theorem}[section]
\newtheorem{theorem}[thm]{Theorem}

\newtheorem{lemma}[thm]{Lemma}
\newtheorem{prop}[thm]{Proposition}

\newtheorem{defn}[thm]{Definition}
\newtheorem{definition}[thm]{Definition}
\newtheorem{questions}{Question}
\newtheorem{conjecture}[thm]{Conjecture}

\newtheorem{rem1}[thm]{Remark}
\newenvironment{remark}{\begin{rem1}\em}{\end{rem1}}

\newcommand{\cnums} {{\mathbf C}}          
\newcommand{\rnums} {{\mathbf R}}		
\newcommand{\znums} {{\mathbf Z}}		

\newcommand{\Hom}{\operatorname{Hom}}
\newcommand{\Tor}{\operatorname{Tor}}
\newcommand{\Ker}{\operatorname{Ker}}
\newcommand{\End}{\operatorname{End}}

\newcommand{\Coker}{\operatorname{Coker}}

\newcommand{\ie}{{\em i.e. }}

\newcommand{\til}[1]{\widetilde{#1}}

\renewcommand{\P}{\mathbf{P}}

\newcommand{\E}{\mathcal{E}}
\renewcommand{\O}{\mathcal{O}}

\renewcommand{\hat}[1]{\widehat{#1}}

\newcommand{\Sym}{\operatorname{Sym}}
\newcommand{\Hilb}{\operatorname{Hilb}}
\newcommand{\sheafext}{\operatorname{\mathcal{E}\!\mathit{xt}}}
\newcommand{\sheafhom}{\operatorname{\mathcal{H}\!\mathit{om}}}
\renewcommand{\dim}{\operatorname{dim}}

\begin{document}

\begin{abstract}
We study the moduli space $M_{G} (A)$ of flat $G$-bundles on an Abelian
surface $A$, where $G$ is a compact, simple, simply connected, connected
Lie group. Equivalently, $M_{G} (A)$ is the (coarse) moduli space of
$s$-equivalence classes of holomorphic semi-stable $G^{\cnums }$-bundles
with trivial Chern classes.

$M_{G} (A)$ has the structure of a hyperk\"ahler orbifold. We show that
when $G$ is $Sp(n)$ or $SU (n)$, $M_{G} (A)$ has a natural hyperk\"ahler
desingularization which we exhibit as a moduli space of $G^{\cnums
}$-bundles with an altered stability condition. In this way, we obtain the
two known families of hyperk\"ahler manifolds, the Hilbert scheme of points
on a $K3$ surface and the generalized Kummer varieties. We show that
for $G$ not $Sp (n)$ or $SU (n)$, the moduli space $M_{G} (A)$ does
\emph{not} admit a hyperk\"ahler resolution.

\sloppy{Inspired by the physicists Vafa and Zaslow, Batyrev and Dais define
``stringy Hodge numbers'' for certain orbifolds. These numbers are
conjectured to agree with the Hodge numbers of a crepant resolution (when
it exists). We compute the stringy Hodge numbers of $M_{SU (n)} (A)$ and
$M_{Sp (n)} (A)$ and verify the conjecture in these cases.}
\end{abstract}


\maketitle 
\markboth{$G$-bundles and hyperk\"ahler manifolds}{$G$-bundles and hyperk\"ahler manifolds}   
\renewcommand{\sectionmark}[1]{}


\section{Results and the motivating examples.}

Recent advances in certain string theories have inspired a resurgence of
interest in the moduli space of $G$-bundles on elliptic curves
\cite{Donagi1} \cite{Donagi2} \cite{Fr-Mo-Wi1} \cite{Fr-Mo-Wi2}
\cite{Fr-Mo-Wi3} \cite{Laszlo}. In these studies, care has been taken to
develop methods that apply to arbitrary $G$ and that are well suited to
families of elliptic curves---the situation of physical interest is
principal bundles on elliptic fibrations with structure group contained in
$E_{8}\times E_{8}$.

In this paper we study flat $G$-bundles on an Abelian surface $A$. We are
primarily interested in the geometry of $M_{G} (A)$, the coarse moduli
space, and so we will not address the existence of a universal family or
the variation of $M_{G} (A)$ in families. This affords us the opportunity
to keep the discussion of $M_{G} (A)$ very concrete and elementary; we have
strived to give the paper some expository value in addition to reporting
our findings.

Before we begin, we summarize our results in the following theorem, deferring
definitions, explanations, and details to the rest of the paper.

\pagebreak

\begin{thm}\label{thm: summary of results}
Let $G$ be a compact, simple, simply connected Lie group. Let $M_{G} (A)$
be the moduli space of flat $G$ bundles on an Abelian surface $A$. Then
\begin{enumerate}
\item $M_{G} (A)$ has a hyperk\"ahler resolution if and only if $G$ is $SU
(n)$ or $Sp (n)$ (Theorem~\ref{thm: M(A) does not have a crepant resolution
for G bad});
\item In these cases, the resolution is realized as a certain moduli space
of $G$-bundles, namely the moduli space of Mukai-stable (see
Definition~\ref{defn: Mukai stability}) $G^{\cnums }$-bundles
(Theorems~\ref{thm: tilMsu is the Mukai stable moduli space} and \ref{thm:
mukai stable Sp bundles is hyperkahler and birat to HilbX}).
\item The stringy Hodge numbers of $M_{SU (n)} (A)$ and $M_{Sp (n)} (A)$
coincide with the ordinary Hodge numbers of their corresponding
hyperk\"ahler resolutions (Theorems \ref{thm: hst(X(n))=h(X[n]) (Gottsche's
computation)} and \ref{thm: hst(Msp)=h(tilMsp)}).
\end{enumerate}
\end{thm}

\subsection{Notation} Fix $A$ to be a principally polarized Abelian surface
(we do this primarily for convenience---our results can be adapted to any
complex torus of dimension 2 without much trouble) and fix $E$ to be an
elliptic curve. We choose origins $p_{0}\in A$ and $p_{0}\in E$. We will
freely identify $A$ and $E$ with their duals $A^{\vee }$ and $E^{\vee }$
(using the polarizations). Let $G$ be a compact, simple, simply connected,
connected Lie group (e.g. $SU (n)$ or $Sp (n)$). Let $G^{\cnums }$ be the
complexification of $G$ (e.g. $SL (n,\cnums )$ or $Sp (n,\cnums )$), let
$r$ be the rank of $G$, and let $W$ be its Weyl group.

\subsection{Questions.} Let $M_{G} (A)$ (respectively $M_{G} (E)$) denote
the moduli space of flat $G$ connections on $A$ (respectively
$E$). Equivalently, $M_{G} (A)$ (respectively $M_{G} (E)$) is the (coarse)
moduli space of $s$-equivalence classes of semi-stable holomorphic
$G^{\cnums }$-bundles on $A$ (respectively $E$) with trivial Chern classes.

In contrast to $M_{G} (E)$, $M_{G} (A)$ is not in general connected. We
denote by $M_{G}^{0} (A)$ the component containing the trivial
connection. This component can be described as a quotient of an $r$-fold
product of $A$ by an action of the Weyl group
\[
M_{G}^{0} (A)\cong A^{r}/W.
\]
The action of $W$ preserves the natural holomorphic symplectic form on
$A^{r}$ and so $M_{G}^{0} (A)$ has a holomorphic symplectic form on the
open dense locus of $W$ orbits with trivial stabilizer (see Sections
\ref{sec: G-bundles on E and A} and \ref{sec: hyperk mfds, symp var.s and
resolutions} for the details of these assertions).

 The questions that
motivated this work are:

\begin{questions}\label{question: does there exist a resolution?}
Does $M_{G}^{0} (A)$ have a smooth resolution $\til{M}_{G}^{0} (A)$
to which the holomorphic symplectic form extends? Such a resolution would
admit a hyperk\"ahler metric.
\end{questions}
\begin{questions}\label{question: is the resolution a moduli space?}
If $\til{M}_{G}^{0} (A)$ exists, can it be realized as a moduli space
for some moduli problem related to $G^{\cnums }$-bundles on $A$?
\end{questions}
\begin{questions}\label{question: how are the geometric structures related?}
How are the  Hodge numbers of the
desingularization $\Tilde{M}_{G}^{0} (A)$ (if it exists) encoded in the
action of $W$ on $A^{r}$?
\end{questions}

The answer to the first two questions is ``yes'' in the case when $G$ is
$Sp (n)$ or $SU (n)$. In these cases, $M^{0}_{G} (A)=M_{G} (A)$ and the
hyperk\"ahler manifolds obtained are exactly the two known families of
irreducible hyperk\"ahler manifolds. $\til{M}_{Sp (n)} (A)$ is
$\operatorname{Hilb}^{n} (X)$, the Hilbert scheme of $n$ points on $X$, the
Kummer $K3$ surface associated to $A$. $\til{M}_{SU (n)} (A)$ is
$KA_{n-1}$, the so called generalized Kummer variety which is the fiber of
the map $\operatorname{Hilb}^{n} (A)\to A$ given by summing the points
using the group law of $A$. We realize these resolutions as the moduli
spaces of ``Mukai-stable'' $G^{\cnums }$-bundles (see Definition \ref{defn:
Mukai stability}). 

A framework for answering the third question is nicely provided by the
``stringy Hodge numbers''. These can be computed purely from the group
theory that defines the action of $W$ on $A^{r}$. In the case of $Sp (n)$
and $SU (n)$ we prove that they give exactly the Hodge numbers of the
resolution.

\subsection{The case of $SU (n)$} Since the case of $SU (n)$ was one of the
motivating examples, we describe it in more detail. The rank of $SU (n)$ is
$n-1$ and we have
\[
A^{n-1}\hookrightarrow A^{n}
\]
as the set of points $(x_{1},\dots ,x_{n})$ with $\sum x_{i}=0$. The Weyl
group $W$ is the symmetric group $S_{n}$ and its action on $A^{n-1}$ is the
restriction of the natural action on $A^{n}$. The identification
\[
M_{SU (n)} (A)\cong A^{n-1}/S_{n}
\]
is easy to understand in concrete terms. Points of $M_{SU (n)} (A)$
naturally correspond to $s$-equivalence classes of holomorphic semi-stable
$SL_{n}$ bundles with trivial Chern classes, that is bundles
$\E\to A$ with $c_{1} (\E)=c_{2} (\E)=0$ and an
isomorphism $\operatorname{det}\E\cong \mathcal{O}_{A}$. In this
case, every semi-stable bundle is strictly semi-stable and $\E$
can be decomposed (up to $s$-equivalence) into a sum of flat line bundles:
\[
\E\cong L_{x_{1}}\oplus \dots \oplus L_{x_{n}}
\]
where $L_{x}$ is the line bundle corresponding to $x\in A\cong
\operatorname{Pic}^{0}A$. This decomposition is unique up to $s$-equivalence
and reordering the factors. The condition that
$\operatorname{det}\E\cong \mathcal{O}$ imposes the condition
$\sum x_{i}=0$.

The singular points of $M_{SU (n)} (A)$ occur on the $S_{n}$-orbits with a
non-trivial stabilizer. This occurs when two or more of the line bundles in
the above description coinciding. When this happens, $s$-equivalence is
rather brutal. It identifies many non-isomorphic bundles to a single moduli
point. To illustrate, consider $SL (2,\cnums )$ bundles on $A$. The moduli
space is
\[
M_{SU (2)} (A)\cong A/\pm 1
\]
where the orbit $\{x,-x \}$ corresponds to the bundle $L_{x}\oplus
L_{-x}=L_{x}\oplus L_{x}^{-1}$. The singular points occur for the sixteen
two torsion points of $A$ where $x=-x$. For a two torsion point $\tau $,
the moduli point $\{\tau ,\tau \}\in A/\pm 1$ corresponds to the
$s$-equivalence class of
\[
L_{\tau }\otimes (\mathcal{O}\oplus \mathcal{O}).
\]

For any non-trivial extension
\[
0\to \mathcal{O}\to \E\to \mathcal{O}\to 0
\]
the bundle $L_{\tau }\otimes \E$ is $s$-equivalent to $L_{\tau
}\otimes (\mathcal{O}\oplus \mathcal{O})$. The natural parameter space for
isomorphism classes of non-trivial extensions of $\mathcal{O}$ by
$\mathcal{O}$ is
\[
\P (\operatorname{Ext}^{1} (\mathcal{O},\mathcal{O}))\cong \P (H^{1}
(A,\mathcal{O}))\cong \P ^{1}.
\]
This suggests that if one could find a way to ``destabilize'' $L_{\tau
}\otimes (\mathcal{O}\oplus \mathcal{O})$ and remove it from the moduli
problem, then the corresponding moduli space should replace each of the
sixteen double points of $M_{SU (2)} (A)$ with $\P ^{1}$'s. Of course, if
we blow up $A/\pm 1$ at the sixteen double points, we obtain $X$, the
Kummer $K3$ surface which is a solution to Question \ref{question: does
there exist a resolution?} in this case. The above discussion also suggests
a strategy for constructing $S$ as a moduli space in order to answer
Question \ref{question: is the resolution a moduli space?}. Following a
suggestion of Aaron Bertram (that goes back to ideas of Mukai), we carry
this out for $G$ equal to $SU (n)$ or $Sp (n)$ in Section~\ref{sec:
Realization of desing's as moduli spaces}. In that section we define a new
notion of stability (Mukai stability). The moduli space of Mukai stable
bundles is then related to the Hilbert scheme of points on the dual Abelian
surface via the Fourier-Mukai transform. Functorial properties of the
Fourier-Mukai transform allow us to carefully analyze the condition that a
bundle has a symplectic structure (the $Sp (n)$ case) where many subtleties
occur.

\subsection{The general case} Unfortunately, this program does not succeed
in producing new examples of compact hyperk\"ahler manifolds. We prove that
$M_{G} (A)$ admits \emph{no} hyperk\"ahler resolution (in fact, no crepant
resolution) for $G$ not $SU (n)$ or $Sp (n)$ (Theorem~\ref{thm: M(A) does
not have a crepant resolution for G bad}). This situation has an analogue
for the moduli space $M_{G} (E)$ of bundles on the elliptic curve $E$. In
\cite{Looijenga}, Looijenga proves that $M_{G} (E)$ is a weighted
projective space. The weighted projective space is smooth if and only if
$G$ is $SU (n)$ or $Sp (n)$. As we will explain, the same mechanism that
causes $M_{G} (E)$ to fail smoothness, causes $M_{G} (A)$ to not admit a
crepant resolution. This analogy continues to hold when we replace $E$
with $\cnums $ and $A$ with $\cnums ^{2}$: Chevalley's theorem asserts that
$\cnums ^{r}/W$ is always smooth; a recently announced result of
Bezrukavnikov-Ginzburg claims that $\cnums ^{2r}/W$ always admits a
holomorphic symplectic resolution. Thus the failure of $M_{G} (A)$ ($G\neq
SU (n)$ or $Sp (n)$) to admit a holomorphic symplectic resolution has to do
with global properties of $A$ (like torsion points). We discuss this
analogy and these results further in Section~\ref{sec: hyperk mfds, symp
var.s and resolutions}.

\subsection{Stringy Hodge numbers} When a Calabi-Yau manifold $X$ is acted
on by a finite group $H$ preserving the holomorphic volume form, Batyrev
and Dais (based on ideas of the physicists Vafa
\cite{Vafa-orbifold-numbers} and Zaslow \cite{Zaslow}) define ``stringy
Hodge numbers'' $h_{st}^{p,q} (X,H)$ \cite{Batyrev-Dais}. In particular, if
$X$ is holomorphic symplectic (e.g. $A^{r}$) and the action of $H$
preserves the symplectic form (e.g. $W$ acting on $A^{r}$), then the
numbers $h_{st}^{p,q} (X,H)$ are well defined. The stringy Hodge numbers
are conjectured to coincide with the ordinary Hodge numbers of a crepant
resolution of $X/H$, if it exists (see Conjecture \ref{conj: string hodge
nums are hodge nums of resolution}). This conjecture is part of the
generalized McKay correspondence.

The situations where this conjecture has been tested are somewhat
limited. It has been verified for $\dim X\leq 3$, and for $H$
Abelian. Since the (ordinary) Hodge numbers of the resolutions of $M_{SU
(n)} (A)$ and $M_{Sp (n)} (A)$ are known, the pairs $(A^{r},W)$ provide
higher dimensional examples with non-Abelian group actions where the
conjecture can be tested. This was done for $G=SU (n)$ by G\"ottsche (see
Theorem~\ref{thm: hst(X(n))=h(X[n]) (Gottsche's computation)}); we verify
the conjecture for $Sp (n)$ (Theorem~\ref{thm: hst(Msp)=h(tilMsp)}). To our
knowledge, there are no other higher dimensional, non-Abelian examples
where this conjecture has been verified.


\subsection{Acknowledgements} It is a pleasure to thank Scot Adams, Aaron
Bertram, Tony Pantev, Michael Thaddeus, and Al Vitter for helpful
conversations.

J. Bryan is supported by an Alfred P. Sloan Research Fellowship and NSF
grant DMS-9802612.  R. Donagi is supported in part by an NSF grant
DMS-9802456 as well as a University of Pennsylvania Research Foundation
Grant; he was partially supported by a grant from the Emmy Noether
Institute and the Minerva foundation of Germany, as well as by the Hebrew
University of Jerusalem and RIMS, Kyoto during visits to those
institutions. C. Leung is supported by NSF grant DMS-9803616.

\section{$G$-bundles on elliptic curves and Abelian surfaces.}
\label{sec: G-bundles on E and A}

Much has been written recently concerning flat $G$-bundles/holomorphic
$G^{\cnums }$-bundles on elliptic curves (for example \cite{Donagi1}
\cite{Donagi2} \cite{Fr-Mo-Wi1} \cite{Fr-Mo-Wi2} \cite{Fr-Mo-Wi3}
\cite{Laszlo}). In this section we follow a ``standard'' approach to the
construction of $M_{G} (E)$ and we develop the theory for Abelian surfaces
in parallel. Since we are mainly interested in the geometry of the coarse
moduli space $M_{G} (A)$, we take an elementary approach to its
construction and ignore the issues of the existence of a universal bundle
and the variation of $M_{G} (E)$ in a family.

\begin{defn}\label{def: MG (X) as flat bundles}
Let $M_{G} (X)$ denote the moduli space of flat $G$-bundles on a path
connected space $X$. It is given by
\[
M_{G} (X)=\Hom (\pi _{1} (X),G)/G
\]
where $G$ acts on a representation by conjugation.
\end{defn}

When $X$ is K\"ahler, there is a correspondence between flat $G$-bundles
and certain holomorphic $G^{\cnums }$-bundles. In the case of $E$ and $A$
it is a special case of the famous theorems of Narasimhan-Seshadri and
Donaldson (generalized by Uhlenbeck and Yau \cite{Uhlenbeck-Yau}):
\begin{thm}[Narasimhan-Seshadri, Donaldson]\label{thm: Nar-Shesadri,Donaldson}
$M_{G} (E)$ (respectively $M_{G} (A)$) is isomorphic to the coarse moduli
space of $s$-equivalence classes of semi-stable holomorphic $G^{\cnums
}$-bundles on $E$ (respectively $A$) with vanishing Chern classes. In
particular, $M_{G} (E)$ and $M_{G} (A)$ are projective varieties.
\end{thm}

See \cite{Huybrechts-Lehn-book} or Section \ref{sec: Realization of
desing's as moduli spaces} for the definitions of semi-stable and
$s$-equivalence. For the most part we will work with the topological
description of these moduli spaces, but we will identify and use the
holomorphic structure coming from the above theorem.

\subsection{Reduction to a finite quotient} Choose a maximal torus
$T\subset G$. By a classical result of Borel \cite{Borel}, any pair of
commuting elements in a compact, simply connected Lie group lie in the same
maximal torus and thus can be simultaneously conjugated to the fixed torus
$T$. Noting that $W=N (T)/T$ is the normalizer of $T$ quotiented by $T$, we
have
\begin{align*}
M_{G} (E)&\cong \Hom (\pi _{1} (E),G)/G\\
&\cong \Hom (\pi _{1} (E),T)/W\\
&\cong (T\times T)/W.
\end{align*}

In general, three or more commuting elements \emph{do not} all lie in the
same maximal torus (although it is true for $SU (n)$ and $Sp (n)$), so the
above analysis for $M_{G} (A)$ does not apply. However, the condition that
commuting elements all lie in the same maximal torus is both open and
closed in $M_{G} (A)$ so if we restrict our attention to the connected
component containing the trivial connection, the above argument will apply.
\begin{defn}\label{def: M0 (A)}
Let $M_{G}^{0} (A)\subset M_{G} (A)$ be the connected component containing
the trivial connection.
\end{defn}
By the previous argument, we then have
\begin{align*}
M^{0}_{G} (A)&\cong \Hom (\pi _{1} (A),T)/W\\
&\cong (T\times T\times T\times T) /W.
\end{align*}

The above description does not make the complex structure of $M_{G}^{0}
(A)$ apparent. To do this we define the \emph{coroot lattice} $\Lambda $ by
the kernel of the exponential map to $T$:
\[
0\to \Lambda \to \mathfrak{t}\to T\to 0.
\]

An element $\pi _{1} (A)\to T$ of $\Hom (\pi _{1} (A),T)$ is dual to a
homomorphism
\[
\Hom (T,S^{1})\to \Hom (\pi _{1} (A),S^{1})\cong A^{\vee}\cong A.
\]
The first group is just $\Lambda ^{\vee }$ and so the above homomorphism is an element of $\Lambda \otimes A$. In this way we have a natural isomorphism
\[
\Hom (\pi _{1} (A),T)\cong \Lambda \otimes A.
\]
The action of $W$ on $\Lambda $ induces an action on $\Lambda \otimes A$
and the complex structure of $A$ induces a holomorphic structure on the
quotient. The same discussion applies to $E$ and so we have
\begin{align}
M_{G} (E)&\cong (\Lambda \otimes E) /W\\
M_{G} ^{0} (A)&\cong (\Lambda \otimes A) /W. \nonumber 
\end{align}
Although we will not prove it, this holomorphic structure is the same as the
one determined by Theorem~\ref{thm: Nar-Shesadri,Donaldson}. Since $\Lambda
$ is a rank $r$ lattice, we may choose a $\znums $-basis and write $\Lambda
\otimes A\cong A^{r}$ as we did in the first section.

\subsection{An example of an unusual commuting triple}Before we continue
our study of $M_{G}^{0} (A)$, we give an example (adapted from a talk of
Witten) showing that there are commuting elements in a simply connected Lie
group that do not all lie in the same maximal torus $T$.

\begin{prop}\label{prop: 3 commuting elts---cant conjugate to T}
Consider the following commuting matrices in $SO (8)$:
\begin{align*}
a=&\operatorname{Diag} (+1,-1,+1,-1,+1,-1,+1,-1)\\
b=&\operatorname{Diag} (+1,+1,-1,-1,+1,+1,-1,-1)\\
c=&\operatorname{Diag} (+1,+1,+1,+1,-1,-1,-1,-1).
\end{align*}
Choose lifts $\til{a},\til{b},\til{c}\in Spin (8)$. Then $\til{a}$,
$\til{b}$, and $\til{c}$ are mutually commuting elements of $Spin (8)$ that
do not all lie in a single maximal torus.
\end{prop}
\begin{remark}\label{rem: commuting elts of SO(3) not in on torus}
In a non-simply connected group, it is easy to find even just two commuting
elements that do not lie in a single maximal torus. For example,
$\operatorname{Diag} (-1,-1,+1)$ and $\operatorname{Diag} (+1,-1,-1)$ are a
commuting pair of $SO (3)$ matrices that are in different maximal tori
(they have different axis of rotation), however any lifts of these elements
to the simply connected cover $SU (2)$ will not commute---their commutator
is $-Id$. For an extensive study of commuting pairs and triples see
\cite{Borel-Friedman-Morgan} and also \cite{Schweigert} or
\cite{Kac-Smilga}.
\end{remark}

\textsc{Proof of Proposition~\ref{prop: 3 commuting elts---cant conjugate
to T}:} We first show that $a$, $b$, and $c$ do not lie in the same maximal
torus in $SO (8)$. We then show that the lifts $\til{a}$, $\til{b}$, and
$\til{c}$ mutually commute. The result will then follow since if $\til{a}$,
$\til{b}$, and $\til{c}$ were contained in the same maximal torus in $Spin
(8)$, then $a$, $b$, and $c$ would be contained in the image torus in $SO
(8)$.

 Let $T^{3}$ be the three torus. The elements $a$, $b$, and $c$ determine a
representation $\pi _{1} (T^{3})\to SO (8)$. \ie a flat $SO (8)$
connection. If $a$, $b$, and $c$ were contained in the same maximal torus,
then they could be simultaneously conjugated to $T$, and the associated
flat bundle would correspond to a moduli point in $(T\times T\times
T)/W$. This bundle would hence have deformations as a flat bundle. We will
show that it does not have deformations.

Real line bundles with a flat connection are parameterized by $H^{1}
(T^{3},\znums /2)\cong (\znums /2)^{3}$. Let $\{R_{\alpha } \}$ be the
eight flat line bundles corresponding to elements $\alpha $ of $H^{1}
(T^{3},\znums /2)$. The holonomy of the direct sum connection on the rank
eight bundle
\[
E=\bigoplus _{\alpha \in H^{1} (T^{3},\znums /2)}R_{\alpha }
\]
around the generators of $\pi _{1} (T^{3})$ is given by the matrices $a$,
$b$, and $c$. To show that $E$ has no deformations as a flat $SO (8)$
bundle we compute the deformation space
\[
H^{1} (T^{3},\mathfrak{so} (E)).
\]
The bundle $\mathfrak{so} (E)$ can be described as the skew symmetric
endomorphisms:
\begin{align*}
\mathfrak{so} (E)&\subset \End (E)\\
&=\oplus _{\alpha ,\beta }R_{\alpha
}\otimes R_{\beta }^{\vee }\\
&\cong\oplus _{\alpha ,\beta }R_{\alpha }\otimes
R_{\beta } 
\end{align*}
so that, $\mathfrak{so}
(E)=\oplus M_{\alpha \beta }$ where the sum is over unordered pairs
$(\alpha ,\beta )$ with $\alpha \neq \beta $ and $M_{\alpha \beta }$ is by
definition the rank 1 subbundle of $( R_{\alpha }\otimes R_{\beta })\oplus
( R_{\beta }\otimes R_{\alpha })$ with local sections $(s_{\alpha }\otimes
s_{\beta },-s_{\beta }\otimes s_{\alpha })$. Note that $M_{\alpha
\beta }\cong R_{\alpha }\otimes R_{\beta }$ as a flat line bundle and
$R_{\alpha }\otimes R_{\beta }\cong R_{\alpha +\beta }$ is the trivial
bundle if and only if $\alpha = \beta $. Thus $\mathfrak{so} (E)$ is a
sum of flat line bundles with no trivial factors. Our claim will then
follow when we show that $H^{1} (T^{3},R_{\alpha })=0$ if $\alpha \neq 0$.

Viewing $T^{3}$ as $S^{1}\times S^{1}\times S^{1}$ we can decompose $\alpha
$ as $(\alpha _{1},\alpha _{2},\alpha _{3})$ by the Kunneth
theorem. Then
\[
R_{\alpha }\cong \pi _{1}^{*} (R_{\alpha _{1}})\otimes \pi _{2}^{*}
(R_{\alpha _{2}})\otimes \pi _{3}^{*} (R_{\alpha _{3}})
\]
where $\pi _{i}$ is the projection on to the $i$th factor and $R_{\alpha
_{i}}$ is the line bundle corresponding to $\alpha _{i}\in H^{1}
(S^{1},\znums /2)$. We then have (again by the Kunneth theorem)
\begin{align*}
H^{1} (T^{3},R_{\alpha })&\cong H^{1} (S^{1},R_{\alpha _{1}})\otimes H^{0} (S^{1},R_{\alpha _{2}})\otimes H^{0} (S^{1},R_{\alpha _{3}})\\
&\oplus H^{0} (S^{1},R_{\alpha _{1}})\otimes H^{1} (S^{1},R_{\alpha _{2}})\otimes H^{0} (S^{1},R_{\alpha _{3}})\\
&\oplus H^{0} (S^{1},R_{\alpha _{1}})\otimes H^{0} (S^{1},R_{\alpha
_{2}})\otimes H^{1} (S^{1},R_{\alpha _{3}}).
\end{align*}
Now $H^{0} (S^{1},R_{\alpha _{i}})=0$ for $R_{\alpha _{i}}$ non-trivial and
$\dim H^{0} (S^{1},R_{\alpha _{i}})=\dim H^{1}
(S^{1},R_{\alpha _{i}})$ by the index theorem and so $H^{1}
(S^{1},R_{\alpha _{i}})=0$ for $R_{\alpha _{i}}$ non-trivial. Thus $H^{1}
(T^{3},R_{\alpha })=0$ unless $\alpha = (\alpha _{1},\alpha _{2},\alpha
_{3})=0$ and so we conclude that $H^{1} (T^{3},\mathfrak{so} (E))=0$.

Finally, the lifts $\til{a}$, $\til{b}$, and $\til{c}$ mutually commute if
and only if the bundle $E$ is spin. We compute $w_{2} (E)$ by the Whitney
product formula:
\[
w (E)=\prod _{\alpha \in H^{1} (T^{3};\znums /2)} (1+\alpha )
\]
and so
\[
w_{2} (E)=\sum _{\alpha ,\beta \in H^{1} (T^{3};\znums /2)}\alpha \cup \beta =0
\]
by the skew-symmetry of the cup product on $H^{1} (T^{3};\znums /2)$. Thus
$E$ is spin and the proposition is proved.  \qed

We leave it as an exercise to the reader to translate the above argument
into a purely algebraic proof.

\subsection{Looijenga's theorem}We return to our study of $M_{G} (E) $ and
$M_{G}^{0} (A)$. The geometry of $M_{G} (E)$ is completely determined by
Looijenga's theorem:
\begin{thm}[Looijenga \cite{Looijenga}]\label{thm: looijengas thm}
$M_{G} (E)\cong (\Lambda \otimes E)/W$ is isomorphic to a weighted
projective space $\P (1,g_{1},\dots ,g_{r})$ where the weights $g_{i}$ are
the coefficients of the highest coroot expressed in terms of the simple
coroots (see Table~1).
\begin{table}\label{table: coefs of highest coroot}
\begin{tabular}{||c|c||}
\hline 
$G$ &	$(g_{1},\dots ,g_{r})$\\ \hline \hline
$SU (n)$&	$(1,\dots ,1)$\\ \hline
$Sp (n)$&	$(1,\dots ,1)$\\ \hline
$Spin (2n)$&	$(1,1,1,2,\dots ,2)$\\ \hline
$Spin (2n+1)$&	$(1,1,2,\dots ,2)$\\ \hline
$G_{2}$&	$(1,2)$\\ \hline
$F_{4}$&	$(1,2,2,3)$\\ \hline
$E_{6}$&	$(1,1,2,2,2,3)$\\ \hline
$E_{7}$&	$(1,2,2,2,3,3,4)$\\ \hline
$E_{8}$&	$(2,2,3,3,4,4,5,6)$ \\ \hline 
\end{tabular}
\vspace{.1in}
\caption{Coefficients of the highest coroots in terms of the simple coroots.}
\end{table}
In particular, $(\Lambda \otimes E)/W$ is a smooth projective
space $\cnums \P ^{r}$ if and only if $G$ is  $SU (n)$ or $Sp (n)$.
\end{thm}

Notice that the theorem fails for $G$ not  simple. For example, if $G=U
(n)$, then $\Lambda_{U (n)} \cong \znums ^{n}$ and $W_{U (n)}\cong S_{n}$
acting on $\znums ^{n}$ by permuting the factors. Thus
\begin{align*}
M_{U (n)} (E)&\cong E^{n}/S_{n}\\
&=\operatorname{Sym}^{n} (E)
\end{align*}
is the $n$th symmetric product of $E$. However, we have an inclusion
$\Lambda _{SU (n)}\subset \Lambda _{U (n)}$ as the rank $n-1$ sublattice of
points $e_{1},\dots ,e_{n}$ with $\sum e_{i}=0$ and the action of $W_{SU
(n)}=W_{U (n)}$ on $\Lambda _{SU (n)}$ is the restriction of the action on
$\Lambda _{U (n)}$. Thus we see that $M_{SU (n)} (E)$ is the fiber over the
origin $p_{0}$ of the sum map
\[
\begin{diagram}
M_{SU (n) } (E)&\rTo&\operatorname{Sym}^{n} (E)\\
&&\dTo>{\operatorname{sum}}\\
&&E
\end{diagram}
\]
By viewing $\operatorname{Sym}^{n} (E)$ as the space of effective degree
$n$ divisors on $E$ and using the canonical isomorphism $E\cong
\operatorname{Pic}^{n}E$, the above sum map is identified with the
Abel-Jacobi map. Then the fiber $M_{SU (n)} (E)$ gets identified with the
linear system $|\mathcal{O} (np_{0})|$ which is indeed a projective space
of dimension $n-1$ as predicted by Looijenga's theorem.

The predicted isomorphism $M_{Sp (n)} (E)\cong \P ^{n}$ arises in a
slightly different way. In this case $\Lambda _{Sp (n)}\cong \znums ^{n}$
and $W_{Sp (n)}$ is a semi-direct product of $ S_{n}$ and $ \{\pm 1\}^{n}$.
$W$ acts on $\Lambda $ by permuting the factors and multiplying each factor
by $\pm 1$ (see the proof of Theorem~\ref{thm: Mspin(7) has a C8/-1 point}
for a detailed discussion of the coroot lattice and $W$ action in this
case). Thus
\begin{align*}
M_{Sp (n)} (E)&\cong E^{n}/ (S_{n}\ltimes \{\pm 1\}^{n})\\
&\cong (E/\pm 1)^{n}/S_{n}\\
&\cong \operatorname{Sym}^{n} (\P ^{1})\\
&\cong \P ^{n}.
\end{align*}

\subsection{Examples: $M_{SU (n)} (A)$ and $M_{Sp (n)} (A)$} We can apply
the previous analysis in the case of the Abelian surface $A$. For these
cases ($Sp (n)$ and $SU (n)$), any collection of commuting matrices is
contained in the same maximal torus. Thus $M^{0}_{Sp (n)} (A)=M_{Sp (n)}
(A)$ and $M_{SU (n)}^{0} (A)=M_{SU (n)} (A)$. As in the elliptic curve
case, we find that $M_{SU (n)} (A)$ is the fiber of the sum map
\[
\begin{diagram}
M_{SU (n)} (A)&\rTo&\operatorname{Sym}^{n} (A)\\
&&\dTo>{\operatorname{sum}}\\
&&A
\end{diagram}
\]
and $M_{Sp (n)} (A)$ is a symmetric product
\[
M_{Sp (n)} (A)\cong \operatorname{Sym}^{n} (A/\pm 1).
\]
Note that unlike for curves, the symmetric product of a surface is
singular. Similarly, while $E/\pm 1\cong \P ^{1}$ is smooth, $A/\pm 1$ is
singular. However, these spaces have natural desingularizations. In general
for a surface $S$, the Hilbert scheme of $n$ points on $S$ together with
the Hilbert-Chow map is a desingularization of $\operatorname{Sym}^{n}
(S)$:
\[
\operatorname{Hilb}^{n} (S)\to \operatorname{Sym}^{n} (S)
\]
(we discuss the Hilbert scheme of points in more detail in Section~\ref{sec:
hyperk mfds, symp var.s and resolutions}).  Likewise, $A/\pm 1$ has a
natural desingularization which is the Kummer $K3$ surface $X$ associated
to $A$. Thus we can construct \emph{ad hoc} desingularizations of $M_{SU
(n)} (A)$ and $M_{Sp (n)} (A)$ as follows. Define $\til{M}_{SU (n)} (A)$ to
be the fiber over $p_{0}$ of the composition $\operatorname{Hilb}^n(A)\to
\operatorname{Sym}^n(A)\to A $ so that we have:
\[
\begin{diagram}
\til{M}_{SU (n)} (A)&\rTo&M_{SU (n)} (A)\\
\dTo&&\dTo\\
\operatorname{Hilb}^{n} (A)&\rTo&\operatorname{Sym}^{n} (A)\\
\dTo&&\dTo>{\operatorname{sum}}\\
A&\rEq&A.
\end{diagram}
\]
Define $\til{M}_{Sp (n)} (A)$ to be $\operatorname{Hilb}^n(X) $ so that we
get the desingularization:
\[
\til{M}_{Sp (n)} (A)=\operatorname{Hilb}^{n} (X)\to \operatorname{Sym}^{n}
(X)\to \operatorname{Sym}^{n} (A/\pm 1)\cong M_{Sp (n)} (A).
\]
In Section~\ref{sec: hyperk mfds, symp var.s and resolutions} we will see
that these \emph{ad hoc} desingularizations are exactly the two known
families (up to deformation) of compact irreducible hyperk\"ahler
manifolds. These give an affirmative answer to Question \ref{question: does
there exist a resolution?} for $Sp (n)$ and $SU (n)$.

In Section~\ref{sec: Realization of desing's as moduli spaces} we will
realize these desingularizations as moduli spaces giving an affirmative
answer to Question \ref{question: is the resolution a moduli space?} for
these cases.

\section{Hyperk\"ahler manifolds, holomorphic symplectic manifolds, and crepant
resolutions.}\label{sec: hyperk mfds, symp var.s and resolutions}

In this section we first give brief expositions of hyperk\"ahler manifolds
and holomorphic symplectic manifolds. A general source for this material
is \cite{Nakajima} and the references therein. We then use some basic facts
about crepant resolutions to determine which $M_{G} (A)$ have
hyperk\"ahler resolutions. The main result of this section is that
$M_{G} (A)$ does \emph{not} admit a hyperk\"aher resolution unless $G$
is $SU (n)$ or $Sp (n)$.

\subsection{Hyperk\"ahler manifolds.}\label{subsec: hyperkahler mflds}

A $4n$ dimensional Riemannian manifold $(X,g)$ is called
\emph{hyperk\"ahler} if the holonomy group of the Levi-Civita connection is
contained in $Sp (n)$. It is called \emph{irreducible hyperk\"ahler} if the
holonomy group is exactly $Sp (n)$. It is well known that up to finite
covers, every hyperk\"aher manifold is a product of irreducible
hyperk\"ahler manifolds and flat tori (e.g.  \cite{Beauville}).  The
hyperk\"ahler condition is equivalent to the existence of a triple of
almost complex structures $(I,J,K)$ such that each is integrable and the
metric $g$ is K\"ahler with respect to any of these structures, and
$(I,J,K)$ satisfy the algebra of the quaternions; that is

\[
\nabla I=\nabla J=\nabla K=0,
\]
and 
\[
I^{2}=J^{2}=K^{2}=IJK=-1.
\]
In fact, there is a whole
2-sphere of K\"ahler structures: for each $(a,b,c)$ with
$a^{2}+b^{2}+c^{2}=1$, the almost complex structure $\lambda =aI+bJ+cK$ is
integrable and $g$ is K\"ahler with respect to $\lambda $. This family of
K\"ahler structures is called the twistor family (c.f. \cite{Br-Le1}).

The holonomy condition imposes very restrictive conditions on the Hodge
theory of a compact hyperk\"ahler manifold $X$. Since $Sp (n)\subset SU
(2n)$, hyperk\"ahler manifolds are Ricci flat and so $h^{0,2n} (X)=1$. In
fact, the whole Hodge diamond is ``mirror symmetric''; that is,
\[
H^{p,q} (X)\cong H^{2n-p,q} (X).
\]
This isomorphism is obtained by wedging a harmonic $(p,q)$ form with a
holomorphic symplectic form (see below) $n-p$ times \cite{Huybrechts}.

Examples of hyperk\"ahler manifolds can be obtained from other
hyperk\"ahler manifolds by a process analogous to symplectic
reduction. Suppose a hyperk\"ahler manifold admits an action of a compact
Lie group $G$ preserving $(g,I,J,K)$, then Hitchin et. al. \cite{HKLR}
introduced the notion of a hyperk\"ahler moment map
\[
\mu :X\to \rnums ^{3}\otimes \mathfrak{g}^{*}
\]
and under suitable conditions, they show that the quotient $\mu ^{-1}
(\zeta )/G$ has a natural induced hyperk\"ahler structure (for example, see
\cite{Kronheimer-ALE}). However, no known, non-trivial examples of this
type are compact unless the original hyperk\"ahler manifold and group are
both infinite dimensional.

\subsection{Holomorphic symplectic manifolds}\label{subsec: compact
hyperkahler and holo'c sympl mnfds} A K\"ahler manifold $X$ of complex
dimension $2n$ is a \emph{holomorphic symplectic manifold} if there exists
a closed, non-degenerate holomorphic 2-form $\sigma \in H^{0} (X,\Omega
^{2}_{X})$. Non-degenerate means that $\sigma ^{n}$ is a non-vanishing
section of $\Omega ^{2n}_{X}=K_{X}$. A holomorphic symplectic manifold is
called \emph{irreducible} if $h^{0} (X,\Omega ^{2}_{X})=1$. The following
is due to Beauville \cite{Beauville}:

\begin{thm}\label{thm: compact mfd is hyperk iff it is holo'c sympl}
A compact manifold $X$ has an (irreducible) hyperk\"ahler metric if and
only if it has a metric such that it is an (irreducible) holomorphic
symplectic manifold.
\end{thm}
\begin{remark}\label{rem: Guan's examples}
If one removes the K\"ahler condition in the definition of holomorphic
symplectic, then this theorem no longer holds. Examples of compact
(non-K\"ahler) complex manifolds with holomorphic symplectic forms and no
hyperk\"ahler structure were constructed by Guan
\cite{GuanII}\cite{GuanIII}.
\end{remark}

\textsc{Sketch of proof of Theorem~\ref{thm: compact mfd is hyperk iff it
is holo'c sympl}:} Suppose that $(X,g)$ is hyperk\"ahler and let $(\omega
_{I},I)$, $(\omega _{J},J)$, and $(\omega _{K},K)$ be the defining K\"ahler
structures. Then with respect to the K\"ahler structure $(\omega _{I},I)$,
it is easy to check that the form $\sigma =\omega _{J}+i\omega _{K}$ is a
holomorphic symplectic form. Conversely, suppose $(X,g)$ is a holomorphic
symplectic manifold with K\"ahler form $\omega $ and holomorphic symplectic
form $\sigma $. Then $\sigma ^{n}$ defines a trivialization of $K_{X}$ and
so by Yau's solution to the Calabi conjecture \cite{Yau}, there is a unique
Ricci flat metric for which $\omega $ is K\"ahler. This gives a reduction
of the holonomy group from $U (2n)$ to $SU (2n)$. Since the Ricci curvature
is zero, the standard B\^ochner argument using the Weitzenb\"ock formula
shows that $\nabla \sigma \equiv 0$. Thus the holonomy group is contained in
$SU (2n)\cap Sp (n,\cnums )=Sp (n)$. One can sharpen this argument to
conclude that the two notions of irreducibility coincide.\qed

\subsection{Examples}\label{subsec: examples of hyperk mfds} Theorem
\ref{thm: compact mfd is hyperk iff it is holo'c sympl} can be used to
construct examples of compact hyperk\"ahler manifolds using the Hilbert
scheme of points on a surface.

Let $\operatorname{Hilb}^{n} (X)$ denote the Hilbert scheme parameterizing
0 dimensional subschemes of length $n$ in a smooth projective surface $X$
(a.k.a. the Hilbert scheme of $n $ points). This turns out to be a smooth
projective variety of dimension $2n$ with many beautiful properties (see
the book by G\"ottsche \cite{Gott-book}). There is a proper morphism (the
Hilbert-Chow morphism) from the Hilbert scheme to the symmetric product
\[
\operatorname{Hilb}
^{n} (X)\to \operatorname{Sym}^{n} (X)
\]
that sends a subscheme $Z\subset X$ to its support (with
multiplicities). Via this map, $\operatorname{Hilb}^{n} (X)$ is a smooth
resolution of $\operatorname{Sym}^{n} (X)$.

The exceptional strata of $\operatorname{Hilb}^{n} (X)$ are in general very
complicated, but over the locus in $\operatorname{Sym}^n(X) $ where no more
than two points coincide, $\operatorname{Hilb}^n(X) $ can be described
explicitly: The Hilbert-Chow morphism is an isomorphism on the locus of
configurations of $n$ distinct points; over configurations with exactly two
points coinciding at $x$ the fiber is a $\cnums \P ^{1}$ parameterizing the
lines in $T_{x}X$. Geometrically, $\operatorname{Hilb}^n(X) $ records in
which direction the two points come together. A local model\footnote{We say
that $A\subset X$ \emph{has a local model} or \emph{is locally modeled on}
$A\subset Y$ if there is an analytic neighborhood of $A$ in $X$ that is
complex analytically isomorphic to a neighborhood of $A$ in $Y$. Note that
$A$ could just be a point and if the subspace or the ambient space is clear
from the context, then we will drop them from the terminology (e.g. ``the
subset $B$ is locally modeled on $Y$'').} for a configuration in
$\operatorname{Sym}^n(X) $ with exactly 2 points coinciding is
\[
\operatorname{Sym}^2(\cnums ^{2})\times \cnums ^{2n-4} 
\]
and the Hilbert-Chow morphism is locally a product:
\[
\operatorname{Hilb}^2(\cnums ^{2})\times \cnums ^{2n-4}\to
\operatorname{Sym}^2(\cnums ^{2})\times \cnums ^{2n-4}.
\]

Now $\operatorname{Sym}^2(\cnums ^{2})= (\cnums ^{2}\times \cnums ^{2})/S_{2}
$ which, after a linear change of variables, is just $(\cnums ^{2}/ \pm
1)\times \cnums ^{2}$. The rational double point in $\cnums ^{2}/\pm 1$ can
be resolved by blowing up. The resulting space is the total space of the
cotangent bundle of $\cnums \P ^{1}$ and the map
\[
T^{*}\cnums \P ^{1}\to \cnums ^{2}/\pm 1
\]
contracts the zero section to the double point. The resolution
$\operatorname{Hilb}^n(X)\to \operatorname{Sym}^n(X) $ near the locus where
2 points coincide is locally modeled on 
\begin{equation}\label{eqn: local model for Hilb2-->Sym2}
T^{*}\cnums \P ^{1}\times \cnums ^{2n-2}\to (\cnums ^{2}/ \pm 1)
\times \cnums ^{2n-2}.
\end{equation}
Note that this local model has a holomorphic symplectic form since for any
complex manifold $M$, $T^{*}M$ has a canonical holomorphic symplectic form
(a fact analogous to the corresponding fact for real manifolds and real
symplectic forms).

This description enabled Fujiki \cite{Fujiki} and Beauville
\cite{Beauville} to construct examples of compact holomorphic symplectic
manifolds (and hence compact hyperk\"ahler manifolds) from Hilbert schemes
of points.

\begin{thm}\label{thm: X sympl surf-->X[n] is symplect}
If $X$ is an algebraic surface that is holomorphic symplectic, then
$\operatorname{Hilb}^{n} (X)$, the Hilbert scheme of $n$ points on $X$, is
a holomorphic symplectic manifold (of complex dimension $2n$).
\end{thm}
\textsc{Sketch of proof:} Recall that $\operatorname{Hilb}^{n} (X)$ is a
smooth resolution of $\operatorname{Sym}^{n} (X)=X^{n}/S_{n}$. If $X$ has a
holomorphic symplectic form, then $X^{n}$ has a natural holomorphic
symplectic form that is invariant under the action of $S_{n}$. Thus
$\operatorname{Sym}^{n} (X)$ has a holomorphic symplectic form on the open
set of $S_{n}$-orbits with trivial stabilizer. The map
$\operatorname{Hilb}^{n} (X)\to \operatorname{Sym}^{n} (X)$ restricts to an
isomorphism on this set, so we get a holomorphic symplectic form on
$\operatorname{Hilb}^{n} (X)$ defined on the complement of the exceptional
set. We need to show that this form extends to a non-degenerate form on all
of $\operatorname{Hilb}^{n} (X)$. This form can be extended to the
complement of the codimension 2 set where 3 or more points come together
using the canonical symplectic form on the local model (Equation~\ref{eqn:
local model for Hilb2-->Sym2}) on this locus. The form then automatically
extends across the codimension 2 strata (by Hartog's theorem) to a form
$\sigma $. The form is non-degenerate since if $\sigma ^{n}$ had a
non-empty zero set, it would have codimension one, but $\sigma $ is
non-degenerate in codimension two by construction.\qed

\begin{remark}\label{rem: could have used douady space}
The restriction to algebraic surfaces is not necessary. The same argument
applies when $X$ is a non-algebraic, holomorphic symplectic surface if we
replace the Hilbert scheme with the corresponding Douady space. We restrict
to algebraic surfaces for convenience only.
\end{remark}

From the Kodaira-Enriques classification of compact complex surfaces, we
know that if a compact algebraic surface $X$ is holomorphic symplectic,
then $X$ must be either a $K3$ or an Abelian surface. If $X$ is a $K3$
surface, then $h^{0} (\operatorname{Hilb}^{n} (X),\Omega ^{2})=1$
(G\"ottsche \cite{Gott}) so $\operatorname{Hilb}^{n} (X)$ is an irreducible
hyperk\"ahler manifold. Since any two $K3$ surfaces are deformation
equivalent, all the examples produced in this way are deformation
equivalent.

For an Abelian surface $A$, $\operatorname{Hilb}^n(A) $ is not
irreducible. However, one can easily see that the holomorphic symplectic
form is non-degenerate on the fibers of the map
$\operatorname{Hilb}^n(A)\to A $ given by the composition of the
Hilbert-Chow map and the sum map $\operatorname{Sym}^n(A)\to A $. Thus the
fibers of $\operatorname{Hilb}^n(A)\to A $, which are, by definition, the
generalized Kummer varieties $KA_{n-1}$, are holomorphic symplectic. One
can also check that $KA_{n-1}$ are irreducible.

Until recently, the only known examples of compact irreducible
hyperk\"ahler manifolds were deformation or birationally equivalent to
$\operatorname{Hilb}^n(X) $ for a $K3$ surface $X$ or $KA_{n}$ for some
Abelian surface $A$. In particular, all the known examples had the same
Betti numbers\footnote{Birationally equivalent hyperk\"ahler manifolds have
the same Betti numbers \cite{Bat}. In fact, it is believed (but not proven)
that birationally equivalent hyperk\"ahler manifolds are actually
deformation equivalent (and hence diffeomorphic). c.f. \cite{Huybrechts}}
as $\operatorname{Hilb}^n(X) $ or $KA_{n}$. However, O'Grady has recently
constructed an isolated example in dimension 10 that does not have the same
Betti numbers as $\operatorname{Hilb}^{5}X$ or $KA_{5}$ \cite{OGrady}.

\subsection{Crepant resolutions}\label{subsec: crepant resolution}

The hyperk\"ahler manifolds $KA_{n}$ and $\operatorname{Hilb}^n(X) $ appear
as resolutions of orbifolds. As we showed in Section~\ref{sec: G-bundles on
E and A}, these orbifolds are $M^{0}_{G} (A)= A^{n}/W$ where $G$ is $SU (n+1)$
and $Sp (n)$ respectively. The orbifolds of the form $A^{n}/W$ are
holomorphic symplectic in the sense that $A^{n} $ has a holomorphic
symplectic form preserved by $W$. The resolution of $M_{G}^{0} (A)$ that
we seek (and have for $SU (n+1)$ and $Sp (n)$) should have a holomorphic
symplectic form that agrees with the holomorphic symplectic form on the
smooth locus of $M_{G}^{0} (A)$. We call such a resolution a holomorphic
symplectic resolution and since $A^{n}/W$ is projective, such a resolution
is a hyperk\"ahler manifold.

\begin{definition}\label{def: hyperkahler resolution}
Let $M$ be a quasi-projective variety non-singular in codimension 1 with a
holomorphic symplectic form defined on the smooth locus of $M$. We say that
a smooth resolution $\til{M}\to M$ is a \emph{holomorphic symplectic
resolution} if $\til{M}$ has a global holomorphic symplectic form that
agrees with the form pulled back from $M$ on the corresponding locus. If
$M$ is projective, we will also call such a resolution a
\emph{hyperk\"ahler resolution}.
\end{definition}

Hyperk\"ahler resolutions are special cases of crepant resolutions:

\begin{defn}\label{defn: crepant resolution}
Let $M$ be a quasi-projective variety, non-singular in codimension one,
with a holomorphic volume form defined on the smooth locus of $M$
(equivalently, $M$ has a trivial canonical class $K_{M}\cong
\mathcal{O}_{M}$). We say a smooth resolution $\til{M}\to M$ is
\emph{crepant} if $\til{M}$ has a global holomorphic volume form that
agrees with the form pulled back from $M$ on the corresponding locus.
\end{defn}
\begin{remark}
We've restricted our definition of crepant to the case where $M$ has
trivial canonical class to emphasize the analogy with holomorphic
symplectic resolutions. In general, an arbitrary proper, birational
morphism $\phi :Y\to X$ has a \emph{discrepancy divisor} $\Delta
=K_{Y}-\phi ^{*}K_{X}$ and $\phi $ is \emph{crepant} if $\Delta =0$.
\end{remark}

Crepant resolutions may not exist in general. Locally, at an isolated
orbifold point, the issue is:
\begin{questions}\label{question: does Cn/G admit a crepant resoltution?}
Let $H\subset SL (n,\cnums )$ be a finite group. When does $\cnums
^{n}/H$ admit a crepant resolution? What can one say about the geometry of
a resolution if it exists?
\end{questions}

These questions are the subject of study of the so-called generalized McKay
correspondence, and to some extent they motivated the definition of stringy
Hodge numbers (see \cite{Batyrev-Dais} \cite{Reid-proc-of-Kinosaki}
\cite{Reid-bourbaki}). We will return to this topic in more detail in
Section~\ref{sec: stringy hodge nums}.

For this section, we quote some existence and non-existence results for
holomorphic symplectic and crepant resolutions.

Consider $\cnums ^{2r}=\cnums ^{2}\otimes \Lambda $ with its $W$
action. The orbifold $\cnums ^{2r}/W$ is the affine analogue of our moduli
space $M_{G}^{0} (A)\cong A^{r}/W$. In a recent announcement of
Bezrukavnikov and Ginzburg \cite{Bezrukavnikov-Ginzburg}, they construct a
holomorphic symplectic resolution of this space.

\begin{thm}[Bezrukavnikov-Ginzburg]\label{thm: Bezrukavnikov-Ginzburg's thm}
$\cnums ^{2r}/W= (\cnums \otimes \Lambda )/W$ admits a holomorphic
symplectic resolution.
\end{thm}

At first glance, this (announced) theorem suggests that $A^{r}/W$ should
also have a holomorphic symplectic (and hence hyperk\"ahler)
resolution. However,\emph{ not all the orbifold points of $A^{r}/W$ have
local models of the type $\cnums ^{2r}/W$}; there are additional
possibilities arising from the presence of torsion in $A$.

The only non-existence result we need is a very simple one that we have
borrowed from the McKay correspondence literature (it is implied by Theorem
5.4 of \cite{Batyrev-Dais} or see \cite{Reid-bourbaki} example 5.4).

\begin{thm}\label{thm: no crepant resolution of C2n/-1 for n>1}
Let $\znums /2=\{\pm 1 \}$ act by $-1$ on all the factors of $\cnums
^{2d}$, $d>1$. Then $\cnums ^{2d}/\pm 1$ does not admit a crepant
resolution.
\end{thm}

Our main result of this section is the following:
\begin{thm}\label{thm: M(A) does not have a crepant resolution for G bad}
Let $G$ be a compact, simple, simply connected Lie group. Let $M_{G} (A)$
be the moduli space of flat $G$ bundles on an Abelian surface $A$. Then
$M_{G}(A)$ admits a crepant resolution if and only if $G$ is $SU (n)$ or
$Sp (n)$; in particular, $M_{G} (A)$ has a hyperk\"ahler resolution if and
only if $G$ is $SU (n)$ or $Sp (n)$.
\end{thm}

We devote the rest of this section to the proof. To prove the theorem as
stated, it obviously suffices to prove it for $M^{0}_{G} (A)$ since when
$G$ is $Sp (n) $ or $SU (n)$, $M_{G}^{0} (A)=M_{G} (A)$.
 
\subsection{The basic examples: $G_{2}$, $B_{3}$, and $D_{4}$.} To prove
Theorem~\ref{thm: M(A) does not have a crepant resolution for G bad}, we
first prove it in the cases when $G$ is $G_{2}$, $Spin (7)$, and $Spin
(8)$, which in Cartan's classification, corresponds to the Dynkin diagrams
$G_{2}$, $B_{3}$, and $D_{4}$. We will later show how the basic examples
can be propagated to every other $G$ not equal to $SU (n)$ or $Sp (n)$.

Theorem~\ref{thm: M(A) does not have a crepant resolution for G bad} for
$G=G_{2}$ follows from Theorem~\ref{thm: no crepant resolution of C2n/-1
for n>1} and the following:

\begin{thm}\label{thm: MG2 has a C4/-1 point}
Let $W$ and $\Lambda $ be the Weyl group and coroot lattice for
$G_{2}$. There exists a point of $(A\otimes \Lambda )/W$ locally modeled on
$\cnums ^{4}/\pm 1$.
\end{thm}
\textsc{Proof:}
$\Lambda $ is the rank two sublattice of $\znums ^{3}$ consisting of those elements summing to zero:
\[
\Lambda =\{(a_{1},a_{2},a_{3})\in \znums ^{3}\text{ :
$a_{1}+a_{2}+a_{3}=0$} \}.
\]
The Weyl group $W$ is the dihedral group of order 12. $W$ is a $\znums
/2$ extension of the symmetric group $S_{3}$ and the $S_{3}$ action on
$\Lambda $ is given by permuting the $a_{i}$'s and the $\znums /2=\{\pm 1
\}$ action is given by $(a_{1},a_{2},a_{3})\mapsto (-a_{1},-a_{2},-a_{3})$.

Choose a triple of distinct 2-torsion points in $A$ that sum to zero; \ie
let $p= (\tau _{1},\tau _{2},\tau _{3})\in A^{3}$ such that $2\tau _{i}=0$,
$\tau _{i}\neq \tau _{j}\neq 0$ for all $i\neq j$, and $\tau _{1}+\tau
_{2}+\tau _{3}=0$. Note that $p\in A\otimes \Lambda $ since
\[
A\otimes \Lambda =\{(x_{1},x_{2},x_{3})\in A^{3}\text{ :
$x_{1}+x_{2}+x_{3}=0$} \}.
\]

Since the $\tau _{i}$ are distinct, no non-trivial permutation fixes $p$;
on the other hand, since $\tau _{i}=-\tau _{i}$, we have that $p=-p$. Thus
the stabilizer of $p$ is $\znums /2=\{\pm 1 \}\subset W$.

Therefore a neighborhood of the image of $p $ in $(A\otimes \Lambda )/W$ is
modeled on $\cnums ^{4}/\pm 1$, where $\pm 1$ actions non-trivially on all
factors.\qed

\begin{remark}
If we replace $A$ with $E$ in the above discussion, we see that in
$M_{G_{2}} (E)= (E\otimes \Lambda) /W$ there is a point modeled on $\cnums
^{2}/\pm 1$. Looijenga's theorem (Theorem~\ref{thm: looijengas thm}) tells
us that $M_{G_{2}} (E)$ is in fact $\cnums \P (1,1,2)$ which has a unique
singular point (modeled on $\cnums ^{2}/\pm 1$). In an elliptic curve $E$,
there are exactly 3 non-zero 2-torsion points and so the choice of the
$\tau _{i}$ is unique (up to permutation). Thus the orbit of $p$ is the
unique singular point in $\cnums \P (1,1,2)$.  In $A$, there are many
choices for the $\tau _{i}$'s and so there are multiple points in
$M_{G_{2}} (A)$ where a crepant resolution does not exist locally.
\end{remark}

The basic examples for $B_{3}$ and $D_{4}$ ($Spin (7)$ and $Spin (8)$) are
variations on the same theme that are slightly more involved. See the
appendix for their construction (Theorems \ref{thm: Mspin(8) has a C8/-1
point} and \ref{thm: Mspin(7) has a C8/-1 point}).

\begin{remark}\label{rem: B3 and G2 examples come from D4}
It is no accident that the examples we have for $D_{4}$, $B_{3}$, and
$G_{2}$ are all very similar. The $D_{4}$ Dynkin diagram has an action of
the symmetric group $S_{3}$ and the ``quotient'' of the $D_{4}$ diagram by
$S_{3}$ ``is'' the $G_{2}$ diagram, while the ``quotient'' of the $D_{4}$
diagram by $\znums /2\subset S_{3}$ ``is'' the $B_{3} $ diagram. What this
really means is that there is an $S_{3} $ action on $\Lambda _{D_{4}}$ so
that $\Lambda _{G_{2}}$ and $\Lambda _{B_{3}}$ are respectively the $S_{3}$
and $\znums /2$ invariant sublattices. In this way, we get an $S_{3}$ action
on $M_{Spin (8)} ^{0}(A)$ so that the $S_{3}$ and $\znums /2$ fixed point sets
give inclusions $M_{G_{2}}^{0} (A)\subset M_{Spin (8)} ^{0}(A)$ and
$M_{Spin (7)}^{0} (A)\subset M_{Spin (8)} ^{0}(A)$. The basic examples for
$G_{2}$ and $B_{3}$ are just the restriction of the basic example for
$D_{4}$ under the above inclusions.
\end{remark}

\subsection{Propagating the basic examples.}  The Dynkin diagrams of
$D_{4}$ and $B_{3}$ (corresponding to $Spin (8)$ and $Spin (7)$
respectively) are 
\vspace{.3in}

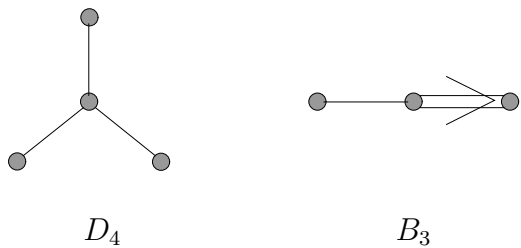
\setlength{\unitlength}{0.00083333in}
\begingroup\makeatletter\ifx\SetFigFont\undefined%
\gdef\SetFigFont#1#2#3#4#5{%
  \reset@font\fontsize{#1}{#2pt}%
  \fontfamily{#3}\fontseries{#4}\fontshape{#5}%
  \selectfont}%
\fi\endgroup%
{\renewcommand{\dashlinestretch}{30}
\begin{picture}(3199,1465)(-1000,-10)
\texture{55888888 88555555 5522a222 a2555555 55888888 88555555 552a2a2a 2a555555 
	55888888 88555555 55a222a2 22555555 55888888 88555555 552a2a2a 2a555555 
	55888888 88555555 5522a222 a2555555 55888888 88555555 552a2a2a 2a555555 
	55888888 88555555 55a222a2 22555555 55888888 88555555 552a2a2a 2a555555 }
\put(511,862){\shade\ellipse{108}{108}}
\put(62,488){\shade\ellipse{108}{108}}
\put(961,487){\shade\ellipse{108}{108}}
\put(513,1389){\shade\ellipse{108}{108}}
\put(1936,862){\shade\ellipse{108}{108}}
\put(2536,862){\shade\ellipse{108}{108}}
\put(3137,863){\shade\ellipse{108}{108}}
\path(100,525)(475,825)
\path(925,525)(550,825)
\path(510,900)(510,1350)
\path(2575,825)(3100,825)
\path(2575,900)(3100,900)
\path(2740,720)(3040,870)(2740,1020)
\path(1975,862)(2500,862)
\put(2425,0){\makebox(0,0)[lb]{\smash{{{\SetFigFont{12}{14.4}{\rmdefault}{\mddefault}{\updefault}$B_3$}}}}}
\put(475,0){\makebox(0,0)[lb]{\smash{{{\SetFigFont{12}{14.4}{\rmdefault}{\mddefault}{\updefault}$D_4$}}}}}
\end{picture}
}

\vspace{.3in}

The $D_{4}$ diagram is a sub-diagram of the diagrams of $E_{6}$, $E_{7}$,
$E_{8}$, and $D_{n}$, $n\geq 4$.  The Dynkin diagram of $B_{3}$ is a
sub-diagram of the diagrams of $F_{4}$ and $B_{n}$, $n\geq 3$.  We can thus
get from the diagrams of $G_{2}$, $B_{3}$, and $D_{4}$ to any other Dynkin
diagram not in the $A_{n}$ or $C_{n}$ series by inclusion.

The following lemma will show that the operation of inclusion allows us to
``propagate'' our basic examples to find points in $M_{G}^{0} (A)$ ($G\neq
Sp (n)$ or $SU (n)$) with local models of the form $(\cnums ^{2l}/\pm
1)\times \cnums ^{2k}$, $l>1$, which then by Theorem~\ref{thm: no crepant
resolution of C2n/-1 for n>1} do not admit crepant resolutions. This will
complete the proof of Theorem~\ref{thm: M(A) does not have a crepant
resolution for G bad}.

\begin{lemma}\label{lem: subdiagrams with bad points cause bad points}
Suppose that $\Lambda \subset \Lambda '$ and $W\subset W'$ are the
inclusions of a coroot lattice and its Weyl group into another coroot
lattice and Weyl group that are induced by an inclusion of a rank $l$
Dynkin diagram into a rank $l+k$ diagram. Let $p\in A\otimes \Lambda $ be a
point and let $W_{p}\subset W$ denote the $W$-stabilizer of $p$. Then there
exists a point $p'\in A\otimes \Lambda '$ such that its $W'$-stabilizer is
isomorphic to $W_{p}$. Moreover, the point $[p']\in ( A\otimes \Lambda ')
/W'$ is locally modeled on $ (\cnums ^{2l}/W_{p})\times\cnums ^{2k}$.
\end{lemma}
\textsc{Proof:} Since $\Lambda '/\Lambda $ is torsion free, the induced map
$A\otimes \Lambda \to A\otimes \Lambda '$ is injective. Via this inclusion,
the $W'$-stabilizer of $p$ (denoted $W'_{p}$) contains the $W$-stabilizer,
\ie $W_{p}\subset W'_{p}$.

Using translation by $p$ and the exponential map, we $W_{p}'$-equivariantly
identify a small neighborhood of $p\in A\otimes \Lambda \subset A\otimes
\Lambda '$ with a small neighborhood of $0\in \cnums ^{2}\otimes \Lambda
\subset \cnums ^{2}\otimes \Lambda '$. Let $N$ be the orthogonal complement
of $\cnums ^{2}\otimes \Lambda $ in $\cnums ^{2}\otimes \Lambda '$. Let $q$
be a small, generic, non-zero element of $N$ which, after exponentiation
and translation, gives us an element $p'\in A\otimes \Lambda '$ lying in a
small neighborhood of $p$. We need to show that $W'_{p'}=W_{p}$;
equivalently we need to show that the $W'_{p}$-stabilizer of $q\in N\subset
\cnums ^{2}\otimes \Lambda '$ (denoted $(W_{p}')_{q}$) is $W_{p}$.

Since $W$ is generated by reflections through planes perpendicular to
vectors in $\Lambda $, elements of $W_{p}\subset W$ fix $N= (\cnums
^{2}\otimes \Lambda )^{\perp }$ and so $(W_{p}')_{q}$ contains $W_{p}$. To
prove the converse, let $g\in (W'_{p})_{q}\subset W_{p}'\subset W'$. Since
$q$ was chosen generically, $g$ must fix all of $N$. We claim that any
element of $W'$ fixing $(\cnums ^{2}\otimes \Lambda )^{\perp }$ must in
fact be an element of $W$. This claim implies that $g\in W'_{p}\cap
W=W_{p}$ and so $(W_{p}')_{q}=W_{p}$ as asserted.

To prove the claim, it is enough to prove the claim with $\cnums ^{2}$
replaced by $\rnums $. In other words, suppose $g\in W'$ acts on
$\mathfrak{t}'=\rnums \otimes \Lambda '$ preserving $\mathfrak{t}^{\perp }=
(\rnums \otimes \Lambda )^{\perp }$; we wish to show that $g\in W$. The set
of Weyl chambers in $\mathfrak{t}'$ (respectively $\mathfrak{t}$) forms a
$W'$-torsor (respectively $W$-torsor). By choosing a fundamental chamber
$C'\subset \mathfrak{t}'$ for $W'$ such that $C=C'\cap \mathfrak{t}$ is a
fundamental chamber for $W$, we get a bijective correspondence between
elements of $W'$ (respectively $W$) and Weyl chambers in $\mathfrak{t}'$
(respectively $\mathfrak{t}$) with the following property.  Those elements
of $W'$ that lie in $W$ correspond to exactly those Weyl chambers in
$\mathfrak{t}'$ whose intersection with $\mathfrak{t}$ is non-trivial. If
$g\in W'$ preserves $\mathfrak{t}^{\perp }$ then it must be an orthogonal
transformation of $\mathfrak{t}$ and so $g (C')\cap \mathfrak{t}=g (C)$ is
a Weyl chamber in $\mathfrak{t}$ and hence $g\in W$ which proves the claim.

Finally, to finish the proof of the Lemma, we observe that via translation
and exponentiation, the decomposition $\cnums ^{2}\otimes \Lambda '= (\cnums
^{2}\otimes \Lambda )\oplus N$ provides the local model whose existence is
asserted by the Lemma. \qed

\section{The stringy Hodge numbers of $M_{SU (n)} (A)$ and $M_{Sp (n)}
(A)$.}  \label{sec: stringy hodge nums}

The Calabi-Yau spaces that appear in the physics of string theory often
have orbifold singularities. Based on mirror symmetry considerations,
physicists have suggested a novel way to extend the definition of Hodge
numbers to Calabi-Yau varieties with orbifold singularities
\cite{Vafa-orbifold-numbers}\cite{Zaslow}.

These so called ``stringy Hodge numbers'' have been extensively studied by
mathematicians recently, especially in the context of the generalized McKay
correspondence (see
\cite{Batyrev-Dais}\cite{Reid-bourbaki}\cite{Reid-proc-of-Kinosaki}). The
stringy Hodge numbers are conjectured to coincide with the ordinary Hodge
numbers of any Crepant resolution, provided it exists.

In this section we compute the stringy Hodge numbers of $M_{Sp (n)} (A)$
and $M_{SU (n)} (A)$ and show that they coincide with the ordinary Hodge
numbers of their resolutions $\til{M}_{Sp (n)} (A)$ and $\til{M}_{SU (n)}
(A)$.

\subsection{Stringy Euler numbers}\label{subsec: stringy Euler nums}
Historically, the stringy Euler number was defined first
\cite{DHVW1}\cite{DHVW2}. Let $X$ be a smooth Calabi-Yau manifold and let
$H$ be a finite group acting on $X$ preserving $K_{X}$

The ordinary Euler number of the quotient $Y=X/H$ can be expressed as
\[
e (X/H)=\frac{1}{|H|}\sum _{g}e (X^{g})
\]
where $X^{g}$ is the fixed point set of an element $g\in H$.

In contrast, the stringy Euler number is \emph{defined} by 
\[
e_{st} (X,H)=\frac{1}{|H|}\sum _{gh=hg}e (X^{g}\cap X^{h})
\]
where the sum is over pairs of commuting elements in $H$.

Note that this sum can be rearranged as a sum over conjugacy classes in
$H$. We let $\{g \}$ denote the conjugacy class represented by an element
$g$. Let $C (g)$ denote the centralizer of $g$. Then $C (g)$ acts on
$X^{g}$ and it is easy to see that $e_{st} (X,H)$ can be rewritten as
\[
e_{st} (X,H)=\sum _{ \{g \}\in \operatorname{Conj} (H)}e (X^{g}/C (g)).
\]

Via localization, it can be shown that $e (X/H)$ is the Euler
characteristic of the $H$-equivariant cohomology of $X$ while $e_{st}
(X,H)$ turns out to be equal to the Euler characteristic of the
$H$-equivariant $K$-theory of $X$ (see \cite{Hirzebruch-Hofer},
\cite{Bryan-Fulman}).

\subsection{Stringy Hodge numbers}\label{subsec: stringy hodge nums} The
stringy Euler number was generalized to Hodge numbers by Zaslow
\cite{Zaslow}. We follow the definition as given by Batyrev and Dais
\cite{Batyrev-Dais}.

Let $X$ be a smooth Calabi-Yau manifold with an action of a finite group
$H$ that preserves the holomorphic volume form. Let $Y=X/H$ be the
quotient. We will define the \emph{stringy Hodge numbers} $h_{st}^{p,q}
(X,H)$, which we will just write as $h_{st}^{p,q} (Y)$ when the orbifold
structure of $Y$ is clear from the context.

\begin{defn}\label{defn: stringy hodge nums}
Let $X$, $H$, and $Y$ be as above. For each $g\in H$, let
\[
X^{g}=X_{1}^{g}\cup \dots \cup X_{r_{g}}^{g}
\]
denote a decomposition of the fixed locus of $g$ into groups of components
which are the orbits of the connected components under the centralizer $C
(g)$ of $g$. For any $x\in X_{i}^{g}$, $g$ acts on $T_{x}X$ with
eigenvalues $e^{2\pi i\alpha _{1}},\dots ,e^{2\pi i\alpha _{n}}$ where we
choose the \emph{weights} $\alpha _{j}$ so that $0\leq \alpha _{j}<1$ and
we define the \emph{Fermion shift number} $F_{i}^{g}$ of the component
$X_{i}^{g}$ as the sum of the corresponding weights; \ie
\[
F_{i}^{g}=\sum _{j=1}^{n}\alpha _{j}.
\]
Note that the $F_{i}^{g}$'s are integers since the action of $g$ preserves
the holomorphic volume form.

We then define the \emph{stringy Hodge numbers} of $Y$ by 
\[
h_{st}^{p,q} (Y)=\sum _{\{g \}\in \operatorname{Conj} (H)}h_{g}^{p,q} (X,H)
\]
where
\[
h_{g}^{p,q} (X,H)=\sum _{i=1}^{r_{g}}h^{p-F_{i}^{g},q-F_{i}^{g}} (X^{g}_{i}/C (g)).
\]
Note that $X_{i}^{g}$ is smooth and has an action of $C (g)$.  The notation
on the right in the above equation is as follows: for any finite group $K$
acting on a smooth complex manifold $V$, let $h^{p,q} (V/K)$ denote the
dimension of the $K$-invariant part of $H^{p,q} (V)$.
\end{defn}

It is easy to see that the stringy Euler characteristic is then given by 
\[
e_{st} (Y)=\sum _{p,q} (-1)^{p+q}h^{p,q}_{st} (Y).
\]

The main conjecture concerning the stringy Hodge numbers is the following.
\begin{conjecture}[Zaslow]\label{conj: string hodge nums are hodge nums of resolution}
If $Z\to X/H$ is any crepant resolution of $X/H$ then
\[
h^{p,q} (Z)=h_{st}^{p,q} (X/H).
\]
\end{conjecture}

\begin{remark}\label{rem: stringy hodge num defn and conj apply to Gorenstien}
The definition of stringy Hodge numbers and the above conjecture are often
both extended to the case where we do not assume that $X$ is Calabi-Yau,
but we merely assume that the singularities of $X/H$ are Gorenstein.
\end{remark}

\subsection{The Fermionic shifts in the holomorphic symplectic case} For a
finite group acting on a holomorphic symplectic manifold preserving the
holomorphic symplectic form, the Fermionic shifts $F_{i}^{g}$ become very
simple.

Recall that the Hodge diamond of a holomorphic symplectic manifold is
mirror symmetric; in other words, if $\dim X=2n$, then the
Hodge diamond is completely symmetric about $(n,n)$ meaning that 
\[
h^{n+p,n+q} (X)=h^{n-p,n-q} (X)=h^{n-p,n+q} (X)=h^{n+p,n-q} (X).
\]

The shifts $F^{g}_{i}$ in the formula for $h_{st}^{p,q} (Y) $ are such that
each individual contribution is completely symmetric (in the above sense)
about $(n,n)$. This can be seen as follows. Since $H$ preserves the
symplectic form, the fixed point sets $X_{i}^{g}$ are each smooth and
holomorphic symplectic. Thus their Hodge diamonds are completely symmetric
about $(n-\frac{1}{2}\operatorname{codim}
(X_{i}^{g}),n-\frac{1}{2}\operatorname{codim} (X_{i}^{g}))$. Thus we just
want to show that
\begin{equation}\label{eqn: shifts are codim(Xg)/2 in the symplectic case}
F_{i}^{g}=\frac{1}{2}\operatorname{codim} (X_{i}^{g}).
\end{equation}
This follows from the fact that $g$ acts on $T_{x}X_{i}^{g} $
symplectically: The eigenvalues of any symplectic transformation come in
pairs $\lambda $, $\lambda ^{-1}$ so the weights $\alpha _{j}$ come in
pairs of the form $(\alpha ,1-\alpha )$ or $(0,0)$. The number of non-zero
eigenvalues is exactly the codimension of $X^{g}_{i}$ and so we have that
\[
F_{i}^{g}=\sum _{j=1}^{2n}\alpha _{j}=\frac{1}{2}\operatorname{codim}
(X_{i}^{g})
\]
as asserted.

\subsection{The computation of $h_{st}^{p,q} (\operatorname{Sym}^n(X))$} In
\cite{Gott-stringy}, G\"ottsche showed that for an algebraic surface $X$,
the stringy Hodge numbers of the symmetric product $\operatorname{Sym}^n(X)
=X^{n}/S_{n}$ coincide with the ordinary Hodge numbers of the resolution
$\operatorname{Hilb}^n(X) $, verifying Conjecture \ref{conj: string hodge
nums are hodge nums of resolution} in this case (c.f. Remark \ref{rem:
stringy hodge num defn and conj apply to Gorenstien}).

In order to state his formula and to facilitate our computations in this
section, we introduce some of G\"ottsche's notation:

\begin{defn}\label{defn: notation for partitions, etc.}
Let $P (n)$ be the set of partitions of $n$. We write $\alpha \in P (n)$ as
$(1^{\alpha _{1}},2^{\alpha _{2}},\dots ,n^{\alpha _{n}})$ so that $\alpha
_{i}$ is the number of $i$'s in the partition. Thus $n=\sum _{i}i\alpha
_{i}$ and we put $|\alpha |=\sum _{i}\alpha _{i}$. We will use the
following shorthand:
\begin{align*}
X^{(n)}&:=\operatorname{Sym}^n(X), \\
X^{[n]}&:=\operatorname{Hilb}^n(X),\\
X^{\alpha }&:=X^{\alpha _{1}}\times \dots \times X^{\alpha _{n}}\text{, and}\\
X^{(\alpha )}&:=X^{(\alpha _{1})}\times \dots \times X^{(\alpha _{n})}, 
\end{align*}
where by convention, $X^{(0)}$ or $X^{0}$ is just a single point. For an
element of the symmetric group $g\in S_{n}$, let $\alpha (g)$ denote its
cycle type and note that $g\mapsto \alpha (g)$ defines a bijection between
conjugacy classes of $S_{n}$ and $P (n)$. Let $X$ be smooth with an action
of a finite group $H$. Define the Hodge and stringy Hodge polynomials by
\begin{align*}
h (X,x,y)&:=\sum _{p,q}h^{p,q} (X)x^{p}y^{q}\\
h_{st} (X/H,x,y)&:=\sum _{p,q}h^{p,q}_{st} (X/H)x^{p}y^{q}.
\end{align*}
Recall that when we use the ordinary Hodge number notation for an orbifold
$h^{p,q} (X/H)$ we mean the dimension of the $H$-invariant part of $H^{p,q}
(X)$. We will then also have the corresponding polynomial:
\[
h (X/H,x,y):=\sum _{p,q} h^{p,q} (X/H)x^{p}y^{q}.
\]
Note that the Hodge polynomial is multiplicative:
\[
h (Y\times Z )=h (Y)h (Z).
\]
\end{defn}

If it happens that the Fermionic shift numbers for all the different
components of $X^{g}$ agree (\ie $F^{g}_{1}=\dots =F^{g}_{r_{g}}$) for all
$g$, then the stringy Hodge polynomial can be written as follows:
\begin{align}\label{eqn: formula for hst poly}
h_{st} (X/H)&=\sum _{p,q}h_{st}^{p,q} (X/H)x^{p}y^{q}\nonumber\\
&=\sum _{p,q}\sum _{\{g \}}h^{p-F^{g},q-F^{g}} (X^{g}/C 
(g))x^{p}y^{q}\nonumber\\
&=\sum _{\{g \}} (xy)^{F^{g}}h (X^{g}/C (g)).
\end{align}

In \cite{Gott-stringy}, G\"ottsche computes the stringy Hodge numbers of
$S^{(n)}$ and $A^{n-1}/S_{n}$ (a.k.a. $M_{SU (n)} (A)$). Having computed
the ordinary Hodge numbers of $S^{[n]}$ and $KA_{n-1}$ in previous works
(\cite{Gott-book}\cite{Gottsche-Soergel}), he then verifies Conjecture
\ref{conj: string hodge nums are hodge nums of resolution} in these cases.

\begin{theorem}[G\"ottsche]\label{thm: hst(X(n))=h(X[n]) (Gottsche's computation)}
For any projective surface $S$ we have
\[
h (S^{[n]},x,y)=h_{st} (S^{(n)},x,y)=\sum _{\alpha \in P (n)}
(xy)^{n-|\alpha |} h (S^{(\alpha )},x,y).
\]
In particular, in the notation of Section~\ref{sec: G-bundles on E and A}, 
\[
h (\til{M}_{U (n)} (A))=h_{st} (M_{U (n)} (A)). 
\]
Moreover, 
\[
h (KA_{n-1})=h_{st} (A^{n-1}/S_{n}),
\]
or, in the notation of Section~\ref{sec: G-bundles on E and A},
\[
h (\til{M}_{SU (n)} (A))=h_{st} (M_{SU (n)} (A)).
\]
\end{theorem}

\subsection{The stringy Hodge number of $M_{Sp (n)} (A)$} G\"ottsche's
theorem verifies Conjecture \ref{conj: string hodge nums are hodge nums of
resolution} for $M_{SU (n)} (A)$ (with its resolution $\til{M}_{SU (n)}
(A)$). We wish to do the same for $M_{Sp (n)} (A)$, \ie we need to compute
the stringy Hodge polynomial of $M_{Sp (n)} (A)$ and compare it to the
ordinary Hodge polynomial of $\til{M}_{Sp (n)} (A)$. The result is the
following:
\begin{thm}\label{thm: hst(Msp)=h(tilMsp)}
Let $M_{Sp (n)} (A)$ and $\til{M}_{Sp (n)} (A)$ be as in Section~\ref{sec:
G-bundles on E and A}. Then
\[
h (\til{M}_{Sp (n)} (A))=h_{st} (M_{Sp (n)} (A)).
\]
In particular, Conjecture \ref{conj: string hodge nums are hodge nums of
resolution} holds for $M_{Sp (n)} (A)$ (with its resolution $\til{M}_{Sp
(n)} (A)$).
\end{thm}

\textsc{Proof:} Recall that $M_{Sp (n)} (A)=A^{n}/S_{n}\ltimes\{\pm 1
\}^{n}$. To compute $h_{st}^{p,q} (M_{Sp (n)} (A))$ we first need to
identify the conjugacy classes of $W=S_{n}\ltimes\{\pm 1 \}^{n}$. An
element of $W$ consists of a permutation along with $n$ signs. We give an
overall sign to each cycle in the permutation by multiplying all of the
signs in the cycle together. In this way, each element determines a
splitting of a partition $\alpha \in P (n)$ into two partitions $\alpha
^{+}+\alpha ^{-}=\alpha $, where $\alpha _{i}^{+}$ and $\alpha _{i}^{-}$ are
the number of cycles of length $i$ of positive and negative type
respectively.  The conjugacy classes of $W$ are in bijective correspondence
with the data $(\alpha ^{+},\alpha ^{-})$ (see Carter \cite{Carter}). For
each $(\alpha ^{+},\alpha ^{-})$, we choose the representative group
element to be such that the positive cycles are first, arranged in order of
increasing length, followed by the negative cycles, also arranged in order
of increasing length. Furthermore, each positive cycle should have all
positive signs and each negative cycle should have exactly one minus sign
at the beginning of each cycle. This determines a unique representative $g
(\alpha ^{\pm })$ of each conjugacy class $(\alpha ^{+},\alpha ^{-})$.

We next determine the fixed point set of $g (\alpha ^{\pm })$ acting on
$A^{n}$. The fixed point locus of each positive cycle fixes is the diagonal
in the product of factors that the cycle acts on. The fixed point locus of
each negative cycle is the 2-torsion points in the diagonal in the product
of the factors that the cycle acts on. Thus the fixed point set of $g
(\alpha ^{\pm })$ is isomorphic to a product of copies of $A$ (one for each
positive cycle) and a product of copies of the set of 2-torsion points of
$A$, denoted $A_{2}$, one for each negative cycle. In other words,
\begin{align*}
(A^{n})^{g (\alpha ^{\pm })}&=\prod _{i=1}^{n}A^{\alpha ^{+}_{i}}\times 
A_{2}^{\alpha _{i}^{-}}\\
&=A^{\alpha ^{+}}\times A_{2}^{\alpha ^{-}}.
\end{align*}

We next need to determine the action of the centralizer $C (g (\alpha ^{\pm
}))$ on the fixed set $(A^{n})^{g (\alpha ^{\pm })}$. Elements that commute
with $g (\alpha ^{\pm })$ permute the cycles of $g (\alpha ^{\pm })$ of the
same length and type. It is easy to see that all such permutations can be
realized (possibly with signs). We only need to understand the signs of the
elements of the centralizer acting on the $A^{\alpha _{i}^{+}} $ factors
since the $A_{2}^{\alpha _{i}^{-}}$ factors are not changed by
multiplication by $-1$ (recall that $A_{2}$ is the group of 2-torsion
points in $A$). Assume then that all the $\alpha _{i}^{\pm }$ are zero
except for a single $\alpha ^{+}_{l}=n/l=m$; that is $g=g (\alpha ^{\pm })$
is a permutation consisting of $m$ cycles all of length $l$ and no minus
signs. We claim that $C (g)/\langle g \rangle$ is then $S_{m}\ltimes \{\pm
1 \}^{m}$ acting on $(A^{n})^{g}\cong A^{m}$ in the standard way. The
$S_{m}\subset C (g)$ is generated by elements that exchange two of the $l$
cycles in $g$; they are a product of $l$ disjoint transpositions with all
positive signs. If each of the transpositions in this product is given a
single minus sign at the beginning of each transposition, then this element
is also in $C (g)$. The effect of the action of this element on
$(A^{n})^{g}\cong A^{m}$ is to permute two of the factors and then multiply
one of them by $-1$. One can directly check that elements of this form
generate $C (g)/\langle g\rangle$.

The above discussion easily generalizes to an arbitrary conjugacy class
$(\alpha ^{+},\alpha ^{-})$. The centralizer $C (g (\alpha ^{\pm }))$ acts
independently on each $A^{\alpha _{i}^{+}}$ and $A_{2}^{\alpha ^{-}_{i}}$
factor, acting by $S_{\alpha _{i}^{-}}$ on $A_{2}^{\alpha _{i}^{-}}$ and by
$S_{\alpha _{i}^{+}}\ltimes\{\pm 1 \}^{\alpha _{i}^{+}}$ on $A^{\alpha
_{i}^{+}}$. Thus we have
\begin{align*}
(A^{n})^{g (\alpha ^{\pm })}/C (g (\alpha ^{\pm })) &=\prod _{i=1}^{n}K^{(\alpha _{i}^{+})}\times A_{2}^{(\alpha _{i}^{-})}\\
&=K^{(\alpha ^{+})}\times A_{2}^{(\alpha ^{-})}
\end{align*}
where $K=A/\pm 1$ is the (singular) Kummer surface associated to $A$.

Since the codimension of $(A^{n})^{g (\alpha ^{\pm })}$ is $2n-2|\alpha
^{+}|$ and the action is symplectic, the Fermionic shift number is just
$n-|\alpha ^{+}|$ (see Equation~\ref{eqn: shifts are codim(Xg)/2 in the
symplectic case}).

From Equation~\ref{eqn: formula for hst poly} and the multiplicative
properties of the Hodge polynomial, we then have
\begin{align}\label{eqn: hst(Msp) in terms of h(K) and h(A2)}
h_{st} (M_{Sp (n)} (A))&=\sum _{\alpha ^{\pm }} (xy)^{n-|\alpha ^{+}|}h (K^{(\alpha ^{+})})h (A_{2}^{(\alpha ^{-})})\nonumber \\
&=\sum _{\alpha ^{\pm }} (xy)^{n-|\alpha|}\prod _{i=1}^{n}h (K^{(\alpha _{i}^{+})})h (A_{2}^{(\alpha _{i}^{-})}) (xy)^{\alpha ^{-}_{i}} \nonumber \\
&=\sum _{\alpha \in P (n)} (xy)^{n-|\alpha |}\prod _{i=1}^{n}\left\{\sum
_{\alpha _{i}^{+}+\alpha _{i}^{-}=\alpha _{i}}h (K^{(\alpha _{i}^{+})})h
(A_{2}^{(\alpha _{i}^{-})}) (xy)^{\alpha _{i}^{-}} \right \}.
\end{align}

Let $X$ be the smooth Kummer $K3$ surface associated to $A$, \ie the blowup
of $K$ at the sixteen double points. Note that
\[
h (X)=h (K)+h (A_{2})xy.
\]

To each polynomial $h (x,y)$ with positive, integral coefficients, we can
assign a bigraded vector space where the dimension of the $(p,q)$ graded
piece is the coefficient of $x^{p}y^{q}$ in $h (x,y)$. The $n$th symmetric
power of this vector space is also a bigraded vector space and so gives
rise to a polynomial which we denote $\operatorname{Sym}^n(h (x,y)) $. By
what is essentially a tautology of the definitions, we have
\[
h (X^{(l)})=\operatorname{Sym}^l (h (X)), 
\]
and from well known properties of the symmetric tensor products we also
have
\[
\operatorname{Sym}^{l} (f+g)=\sum
_{l^{+}+l^{-}=l}\operatorname{Sym}^{l^{+}} (f)\operatorname{Sym}^{l^{-}}
(g).
\]

Thus we get
\begin{align*}
h (X^{(\alpha _{i})})&=\Sym ^{\alpha _{i}} (h (X))\\
&=\Sym ^{\alpha _{i}} (h (K)+h (A_{2})xy)\\
&=\sum _{\alpha _{i}^{+}+\alpha _{i}^{-}=\alpha _{i}}\Sym ^{\alpha _{i}^{+}} (h (K))\Sym ^{\alpha _{i}^{-}} (h (A_{2})xy)\\
&=\sum _{\alpha _{i}^{+}+\alpha _{i}^{-}=\alpha _{i}}h (K^{(\alpha
_{i}^{+})})h (A_{2}^{(\alpha _{i}^{-})}) (xy)^{\alpha _{i}^{-}}
\end{align*}
and so substituting into Equation~\ref{eqn: hst(Msp) in terms of h(K) and
h(A2)} and using Theorem~\ref{thm: hst(X(n))=h(X[n]) (Gottsche's
computation)} we get
\begin{align*}
h_{st} (M_{Sp (n)} (A))&=\sum _{\alpha \in P (n)} (xy)^{n-|\alpha |}\prod _{i=1}^{n}h (X^{(\alpha _{i})})\\
&=\sum _{\alpha \in P (n)} (xy)^{n-|\alpha |}h (X^{(\alpha )})\\
&=h (X^{[n]})\\
&=h (\til{M}_{Sp (n)} (A))
\end{align*}
completing the proof of Theorem~\ref{thm: hst(Msp)=h(tilMsp)}.\qed

\subsection{Stringy speculations.}
It turns out that the series
\[
\sum _{n=0}^{\infty }h_{st} (M_{Sp (n)} (A),x,y)q^{n}
\]
and 
\[
\sum _{n=0}^{\infty }h_{st} (M_{SU (n)} (A),x,y)q^{n}
\]
have some interesting arithmetic properties. Certain series obtained by
setting $x$ and $y$ to special values can be expressed in terms of modular
and quasi-modular forms. For example, by setting $(x,y)= (-1,1)$ one
obtains the generating series for the signature, while setting $(x,y)=
(-1,-1)$ one obtains the generating series for the Euler characteristic.
G\"ottsche gives expressions for these series as the Fourier expansions of
certain quasi-modular forms (\cite{Gott-book} pages 37--39, 51--53, and
57).

Regardless of the existence of crepant resolutions, the stringy Hodge
numbers $h_{st}^{p,q} (M_{G} (A))$ are well defined for any $G$. It would
be interesting to compute the generating series for the $B_{n} $ and
$D_{n}$ series (\ie Spin(odd) and Spin(even)) and determine if they also
have nice expressions. This calculation is
straight forward, although rather complicated. We conjecture that the
generating series for the stringy signature and stringy Euler
characteristic have closed expressions in terms of quasi-modular forms.

\section{Realization of the desingularizations as moduli spaces}\label{sec: Realization of desing's as moduli spaces}

The previous desingularizations of $M_{SU (n)} (A)$ and $M_{Sp (n)} (A)$
defined in Section~\ref{sec: G-bundles on E and A} were \emph{ad hoc}. In
this section, we expand on a suggestion of Aaron Bertram (that also is
implicitly contained in the work of Mukai), to realize these
desingularizations as moduli spaces of holomorphic $G^{\cnums }$-bundles
satisfying a stability condition that we call Mukai-stability. This will
give an affirmative answer to Question \ref{question: is the resolution a
moduli space?} (for the $SU (n)$ and $Sp (n)$ cases).

The notion of Mukai-stability refines ordinary semi-stability in the sense
that Mukai stable bundles are semi-stable. Furthermore, the generic
semi-stable bundle is also Mukai stable. Consequently, the map from the
moduli space of Mukai stable bundles to the moduli space of semi-stable
bundles is generically an isomorphism. However, non-generically, the
stability notions differ which results in a different equivalence relation
for the corresponding moduli problems. In this way, the moduli space for
Mukai stable bundles actually becomes a resolution of singularities for the
ordinary moduli space.

To define Mukai-stability we use a ``framing'' of the fiber over $p_{0}$ of
the holomorphic vector bundle associated to a principal $G^{\cnums }$
bundle. That is, if $\E$ is a holomorphic bundle, a
\emph{framing}\footnote{We use the terminology ``framing'' following
Huybrechts and Lehn \cite{Huybrechts-Lehn} who more generally define a
``framed module'' as a pair $(\E ,\tau )$ where $\tau :\E \to \mathcal{F}$
is a sheaf map to some fixed sheaf $\mathcal{F}$.} is a surjective sheaf
map $\tau :\E\to \mathcal{O}_{p_{0}}$ where $\mathcal{O}_{p_{0}}$ is the
one dimensional skyscraper sheaf at the origin $p_{0}$. We say that the
bundle $\E$ is \emph{Mukai-stable} if it is semi-stable and $\Ker (\tau ) $
is simple, \ie the automorphism group of $\Ker (\tau ) $ is $\cnums
^{*}$. It turns out that $\Ker (\tau )$ is independent of $\tau $ (as long
as $\Ker (\tau )$ is simple).

We motivate how this sort of condition arises by considering how the
desingularization $\Hilb ^{n} (A)\to \Sym ^{n} (A)$ occurs in the context
of moduli spaces and then interpreting this in the context of bundles via
the isomorphism $M_{U (n)} (A)\cong \Sym ^{n} (A)$. The main tool is the
Fourier-Mukai transform. A detailed study of this case allows us to
generalize to the cases of $SU (n)$ and $Sp (n)$.

\subsection{Sym, Hilb, and the moduli of sheaves} We begin by discussing
the isomorphism $M_{U (n)} (A)\cong \Sym ^{n} (A)$ from the point of view
of holomorphic bundles rather than flat connections. Recall that we have
been viewing $M_{U (n)} (A)$ as parameterizing flat $U (n)$ connections so
that the isomorphism with $\Sym ^{n} (A)$ was obtained by ``diagonalizing''
the representation $\pi _{1} (A)\to U (n)$ (see Section~\ref{sec: G-bundles
on E and A}). By Donaldson's theorem (Theorem \ref{thm:
Nar-Shesadri,Donaldson}), we may also view $M_{U (n)} (A)$ as
parameterizing $s$-equivalence classes of semi-stable holomorphic $Gl
(n,\cnums )$ bundles, or, equivalently, $M_{U (n)} (A)$ parameterizes
$s$-equivalence classes of semi-stable, rank $n$, holomorphic bundles $\E$
with
\[
ch (\E)= (n,0,0)\in H^{0} (A,\znums )\oplus H^{2} (A,\znums
)\oplus H^{4} (A,\znums ).
\]

More generally, if $v= (n,c,\chi )\in H^{0} (A,\znums )\oplus H^{2}
(A,\znums )\oplus H^{4} (A,\znums )$ is any vector, there exists (Simpson
\cite{Simpson}, c.f. \cite{Huybrechts-Lehn-book}) a (coarse) moduli space,
which we denote $M (v)$, of $s$-equivalence classes of semi-stable,
coherent, pure-dimensional sheaves $\mathcal{F}$ with $ch (\mathcal{F})=
(n,c,\chi )$.

Using this notation for this section, we write $M (n,0,0)$ instead of $M_{U
(n)} (A)$. The isomorphism $M (n,0,0)\cong \Sym ^{n} (A)$ is obtained by
arguing that every bundle $\E$ with $ch (\E)= (n,0,0)$ is $s$-equivalent to
a direct sum of degree 0 line bundles, \ie $L_{x_{1}}\oplus \dots \oplus
L_{x_{n}}$ where $x_{i}\in \operatorname{Pic}^{0} (A)\cong A$. This follows
from basic facts concerning $s$-equivalence and Jordan-H\"older
filtrations: in general, every semi-stable bundle $\E$ admits a filtration
$0=\E _{0}\subset \E _{1}\subset \dots \subset \E _{k}=\E $ where the
sheaves $\E _{i+1}/\E _{i}$ are stable and have constant slope. Let $Gr
(\E)=\oplus _{i}\E _{i+1}/\E _{i}$, then $\mathcal{F}$ is $s$-equivalent to
$\E$ if and only if $Gr (\E)\cong Gr (\mathcal{F})$. In the case at hand,
the semi-stability of $\E $ and the stability and constant slope of the $\E
_{i+1}/\E _{i}$'s imply that all the factors of $Gr (\E)$ must be degree
0 line bundles.

We see that under $s$-equivalence, all the possible different extensions
\[
0\to \E _{i}\to \E _{i+1}\to \E _{i+1}/\E _{i}\to 0
\]
get identified in the same $s$-equivalence class. In order to find a
desingularization, we would like to find a different notion of equivalence
that remembers such information.

To motivate how one should go about doing this, we first consider how the
desingularization
\[
\Hilb ^{n} (A)\to \Sym ^{n} (A)
\]
naturally occurs in the context of moduli spaces of sheaves. $\Hilb ^{n}
(A)$, by definition, parameterizes length $n$, 0-dimensional subschemes of
$A$. Such a subscheme $Z$ is determined by its ideal sheaf
$\mathcal{I}_{Z}$ which can be considered as a point in $M (1,0,n)$ (the
moduli space of rank 1 semi-stable sheaves $\mathcal{F}$ with $c_{1}
(\mathcal{F})=0$ and $c_{2} (\mathcal{F})=-n$). In fact, there is an
isomorphism
\[
\Hilb ^{n} (A)\times A \cong  M (1,0,n)
\]
given by $(Z,x)\mapsto \mathcal{I}_{Z}\otimes L_{x}$. 

Now the subscheme $Z\subset A$ is also determined by its structure sheaf
$\mathcal{O}_{Z}$ which we can regard as a rank 0 sheaf on $A$ where $ch
(\mathcal{O}_{Z})= (0,0,n)$. Thus $\mathcal{O}_{Z}$ determines a point in
$M (0,0,n)$ . However, $M (0,0,n)$ is isomorphic to $\Sym ^{n} (A)$. The
reason is that every sheaf in $M (0,0,n)$ is $s$-equivalent to a sheaf of
the form $\O _{x_{1}}\oplus \dots \oplus \O _{x_{n}}$. For example, if
$Z\subset A$ is a subscheme of length 2 supported at $x\in A$ then $\O
_{Z}$ can be written as a non-trivial extension
\[
0\to \O _{x}\to \O _{Z}\to \O _{x}\to 0
\]
which is $s$-equivalent to the trivial extension $\O _{x}\oplus \O _{x} $. 

Under these isomorphisms, the Hilbert-Chow morphism
\[
\Hilb ^{n} (A)\to \Sym ^{n} (A)
\]
is obtained by sending the moduli point of $\mathcal{I}_{Z}$ to the
moduli point of $\O _{Z}$.

To translate this picture over to rank $n$ bundles we need a correspondence
between sheaves on $A$ and sheaves on $\operatorname{Pic}^{0} (A)$ that
generalizes the tautological correspondence between points in
$\operatorname{Pic}^{0} (A) $ and line bundles on $A$. The Fourier-Mukai
transform provides such a dictionary.

\subsection{The Fourier-Mukai transform} Although we have been identifying
$A$ and $\operatorname{Pic}^{0} (A)$ throughout his paper via the
polarization, for clarity in this section we write $\hat{A}$ for
$\operatorname{Pic}^{0} (A)$. The ideas of this subsection are due to
Mukai; see the papers \cite{Mukai-Duality} and \cite{Mukai-Senai}.

Let $P\to A\times \hat{A}$ denote the normalized Poincare bundle, \ie
$P|_{A\times \{x \}}=L_{x}$ and $P|_{\{p_{0} \}\times \hat{A}}$ is
trivial. Let $\pi :A\times \hat{A}\to A$ and $\hat{\pi }:A\times \hat{A}\to
\hat{A}$ denote the projections.  Define the functors 
\[
S (?)=\hat{\pi }_{*}
(\pi ^{*} (?)\otimes P)
\]
and 
\[
\hat{S} (?)={\pi }_{*} (\hat{\pi } ^{*}
(?)\otimes P).
\]
If $\mathcal{F}$ is a sheaf on $A$ (respectively,
$\hat{A}$), we obtain sheaves $R^{i}S (\mathcal{F})$ (respectively,
$R^{i}\hat{S} (\mathcal{F})$) on $\hat{A}$ (respectively, $A$) via the
right derived functors of $S$ (respectively, $\hat{S}$).

\begin{defn}\label{defn: weak index thm (WIT)}
If $R^{i}S (\mathcal{F})=0$ for all $i$ except some $i_{0}$, then we say
that $\mathcal{F}$ satisfies the \emph{Weak Index Theorem (W.I.T.)} with
index $i (\mathcal{F})=i_{0}$ and we call the sheaf $R^{i_{0}}S
(\mathcal{F})$ the \emph{Fourier-Mukai transform} of $\mathcal{F}$ and
denote it $\hat{\mathcal{F}}$.
\end{defn}

We also have the analogous definition for sheaves on $\hat{A}$ and
furthermore, if $\mathcal{F}$ satisfies W.I.T., then $\hat{\mathcal{F}}$
satisfies W.I.T.  and $\hat{\hat{\mathcal{F}}}= (-1_{A})^{*}\mathcal{F}$.

For example, any degree 0 line bundle $L_{\hat{x}}\to A$ corresponding to
the point $\hat{x}\in \hat{A}$ satisfies W.I.T.  (with $i (L_{\hat{x}})=2$)
and $\hat{L}_{\hat{x}}=\O _{-\hat{x}}$. More generally, every bundle $\E
\to A$ with $ch (E)= (n,0,0)$ satisfies W.I.T. and the Fourier-Mukai
transform induces an isomorphism
\[
\hat{}:M (n,0,0)\to M (0,0,n).
\]
Conversely, the structure sheaf $\O _{Z}$ of a length $n$, zero-dimensional
subscheme $Z\subset A$, satisfies W.I.T. with $i (\O _{Z})=0$ and $\hat{\O
_{Z}}$ is a bundle $\E $ with $ch (\E)= (n,0,0)$.

Note that despite the fact that the moduli spaces $M (n,0,0)$ and
$M(0,0,n)$ (each isomorphic to $\Sym ^{n} (\hat{A})$) both suffer from an
undiscriminating $s$-equivalence, the Fourier-Mukai transform itself does
not lose information. For example, if $\E $ is a non-trivial extension
\[
0\to \O \to \E \to \O \to 0
\]
then $\hat{\E }=\O _{Z}$ where $Z$ is a length 2 subscheme supported at
$p_{0}\in \hat{A}$ and so it is a non-trivial extension of $\hat{\O }=\O
_{p_{0}}$ by $\hat{\O }=\O _{p_{0}}$:
\[
0\to \O _{p_{0}}\to \O _{Z}\to \O _{p_{0}}\to 0.
\]
The $\P ^{1}$'s worth of non-trivial extensions of $\O $ by $\O $
correspond to the $\P ^{1}$'s worth of length 2 subschemes supported at a
point.

It is not quite the case that the Fourier-Mukai transform is an equivalence
of the category of coherent sheaves on $A$ with the category of coherent
sheaves on $\hat{A}$. However, the Fourier-Mukai transform is an
equivalence of the subcategories of sheaves satisfying W.I.T.. More
generally, the derived functor of $S$ defines an equivalence between the
derived categories of coherent sheaves on $A$ and $\hat{A}$:
\[
\mathbf{R}S:D (A)\to D (\hat{A})
\]
(see \cite{Thomas-Derived} for an introduction to the derived category or
\cite{Weibel} for more detail).  The inverse functor is easily determined
because there is an isomorphism of functors, $\mathbf{R}S\circ
\mathbf{R}\hat{S}\cong (-1_{A})^{*}[-2]$ where ``$[-2]$'' denotes ``shift
the complex 2 places to the right''.

\subsection{Translating points to bundles.}  We now examine what happens to
the sheaves $\mathcal{I}_{Z}$ and $\O _{Z}$, their corresponding moduli
spaces, and the exact sequence $0\to \mathcal{I}_{Z}\to \O \to \O _{Z}\to
0$ under the equivalence of categories $\mathbf{R}S$.

We saw in the last subsection that for any length $n$, zero dimensional
subscheme $Z$, $\O _{Z}$ satisfies W.I.T. with $i (\O _{Z})=0$ and
$\hat{\O_{Z}}$ is a bundle $\E $ with $ch (\E)= (n,0,0)$. That is,
$\mathbf{R}S (\O _{Z})$ is represented by the complex of sheaves that in
degree 0 is $\E $ and it is 0 in all other degrees. $\O $ also satisfies
W.I.T. with $i (\O )=2$ and $\hat{\O }=\O _{p_{0}}$; that is, $\mathbf{R}S
(\O )$ is represented by the complex $\O _{p_{0}}[2]$ which is $\O
_{p_{0}}$ in degree 2 and 0 otherwise.

The exact sequence of sheaves
\[
0\to \mathcal{I}_{Z}\to \O \to \O _{Z}\to 0
\]
is an exact triangle in $D (A)$ when we regard the sheaves as complexes
concentrated in degree 0. Since the functor $\mathbf{R}S$ is an equivalence
of categories, it must take exact triangles to exact triangles and so
\[
\mathbf{R}S (\mathcal{I}_{Z})\to \O _{p_{0}}[2]\to \E \to \mathbf{R}S
(\mathcal{I}_{Z})[1]
\]
is an exact triangle in $D (\hat{A})$.

This triangle gives us a long exact sequence in cohomology from which we
see immediately that $\mathbf{R}S (\mathcal{I}_{Z})$ cannot be a sheaf; it
is represented by a complex of sheaves whose cohomology is $\O _{p_{0}}$ in
degree 2, $\E $ in degree 1, and 0 otherwise. The problem is that
$\mathbf{R}S (\O )$ and $\mathbf{R}S (\O _{Z})$ are sheaves, but they are
concentrated in different degrees. To rectify this problem, we employ
another functor that is an equivalence of derived categories.

Let $\Delta (?)=\sheafhom  (?,\O )$ be the
dualizing functor and let
\[
\mathbf{R}\Delta :D (\hat{A})\to D (\hat{A})
\]
be its derived functor (warning: our notation for $\Delta $ differs from
Mukai's, his has the $\mathbf{R}$ built in and has an additional shift of
the index by 2). Since $\O _{A}$ is the dualizing sheaf of $A$, $\mathbf{R}\Delta $
is an anti-equivalence of the category $D (\hat{A})$ to itself. The
composition
\[
\mathbf{R}\Delta \circ \mathbf{R}S=\mathbf{R} (\Delta S):D (A)\to D
(\hat{A})
\]
is thus also an anti-equivalence. The derived dual of $\E$ is simply the
ordinary dual, \ie $\mathbf{R}\Delta (\E )=\E ^{\vee }$. The derived dual
of $\O _{p_{0}}[2]$ is also a sheaf concentrated in degree zero since
\[
R^{i}\Delta (\O _{p_{0}}[2])=\sheafext ^{2+i} (\O
_{p_{0}},\O )\cong 
\begin{cases} \O _{p_{0}}&i=0\\ 0& i\neq 0.\end{cases}
\]

Thus if we apply $\mathbf{R} (\Delta   S)$ to the sequence
$\mathcal{I}_{Z}\to \O \to \O _{Z}$ (which reverses arrows), we get the
exact triangle
\[
\mathbf{R} (\Delta   S) (\mathcal{I}_{Z})\leftarrow \O
_{p_{0}}\leftarrow \E ^{\vee }\leftarrow\mathbf{R} (\Delta   S)
(\mathcal{I}_{Z})[-1].
\]
Since the map $\O \to \O _{Z} $ is non-zero, the map $\E ^{\vee }\to \O
_{p_{0}}$ must be non-zero and hence surjective. Thus we see that
$\mathbf{R} (\Delta   S) (\mathcal{I}_{Z})$ is concentrated in degree
-1 and is the sheaf $\Ker (\E ^{\vee }\to \O _{p_{0}})$.

We have shown that the functor $\mathbf{R} (\Delta   S)$, which is an
(anti-)equivalence of categories, takes the ideal sheaves of
zero-dimensional subschemes of $A$ to sheaves on $\hat{A}$ which are the
kernel of a framing $\tau :\E \to \O _{p_{0}}$. It remains to see which
bundles and framings arise in this way. This can be answered by reversing
the question; when is $\mathbf{R} (\Delta   {S}) (\Ker (\E \to \O
_{p_{0}}))$ a sheaf of the form $\mathcal{I}_{\hat{Z}}$ where
$\mathcal{I}_{\hat{Z}}$ is the ideal sheaf of a 0 dimensional subscheme
$\hat{Z}$ of $\hat{A}$? We can reverse the question in this way because the
functor $\mathbf{R} (\Delta   \hat{S})$ is the inverse of $\mathbf{R}
(\Delta   S)$, \ie $\mathbf{R} (\Delta   \hat{S})\circ \mathbf{R}
(\Delta   S)$ is isomorphic to the identity functor (Theorem~2.2. and
Equation~3.8 of \cite{Mukai-Duality}). The answer to this question is
given by the following theorem of Mukai (c.f.  \cite{Mukai-Senai}
Proposition~2.18 and Corollary~2.19).

\begin{thm}\label{thm: Mukai's thm}
Let $\E$ be a holomorphic bundle on $A$ with $ch (\E)= (n,0,0)$, let $\tau
:\E \to \O _{p_{0}}$ be a surjective sheaf map, and let $\mathcal{F}=\Ker
(\tau )$. Then the following are equivalent:
\begin{enumerate}
\item $\mathbf{R} (\Delta   {S}) (\mathcal{F})$ is a sheaf;
\item $\hat{\E }$, the Fourier-Mukai transform of $\E $, is of the form $\O
_{\hat{Z}}=\O /\mathcal{I}_{\hat{Z}}$ where $\hat{Z}\subset \hat{A}$ is a
length $n$, 0 dimensional subscheme of $\hat{A}$.
\item $\mathcal{F}$ is simple, i.e. $\End (\mathcal{F})=\cnums $;
\end{enumerate}
In this case $\mathbf{R} (\Delta   {S} )
(\mathcal{F})=\mathcal{I}_{\hat{Z}}[-1]$.
\end{thm}

Note that the theorem implies that $\mathcal{F}$ is independent of the
choice of $\tau $ (as long as $\Ker (\tau )$ is simple).

\textsc{Proof:} First we note that $\mathbf{R} (\Delta S) (\E
)=\mathbf{R}\Delta (\mathbf{R}S (\E))=\mathbf{R}\Delta (\hat{\E
}[2])=\hat{\E } $ and $\mathbf{R} (\Delta S) (\O _{p_{0}})=\mathbf{R}\Delta
(\O ) =\O$.  We apply the functor $\mathbf{R} (\Delta S)$ to the exact
sequence $0\to \mathcal{F}\to \E \to \O _{p_{0}}\to 0$ to get the exact
triangle
\[
\mathbf{R} (\Delta S) (\mathcal{F})\leftarrow \hat{\E }\leftarrow \O
\leftarrow \mathbf{R} (\Delta S) (\mathcal{F})[-1]
\]

Suppose that $\mathbf{R} (\Delta S) (\mathcal{F}) $ is a sheaf.  Since
$\hat{\E }$ is a sheaf supported on a finite number of points, $R^{-1}
(\Delta S) (\mathcal{F})=\Ker (\O \to \hat{\E })\neq 0$. Thus $\mathbf{R}
(\Delta S) (\mathcal{F}) $ is a sheaf implies that $R^{0} (\Delta S)
(\mathcal{F})=\Coker (\O \to \hat{\E })=0$. Thus $\O \to \hat{\E}$ is
surjective and since $\hat{\E }$ is supported on points the kernel of $\O
\to \hat{\E}$ must be an ideal sheaf of a zero dimensional subscheme. Thus
(1) implies (2) and $\mathbf{R} (\Delta {S} )
(\mathcal{F})=\mathcal{I}_{\hat{Z}}[-1]$. It follows then that
$\mathcal{F}$ is simple since $\mathcal{I}_{\hat{Z}}$ is a simple sheaf and
$\mathbf{R} (\Delta S)$ is an equivalence of categories. It remain to be
seen that (3) implies (1).

Assume that $\mathcal{F}$ is simple and suppose that $\mathbf{R} (\Delta S)
(\mathcal{F})$ is not a sheaf. Then $R^{0} (\Delta S) (\mathcal{F})=\Coker
(\O \to \hat{\E})\neq 0$. Since $\hat{\E}$ is supported on a finite number
of points, so is $R^{0} (\Delta S) (\mathcal{F})$ and so there must exist
some point $\hat{x}\in \hat{A}$ so that
\[
\Hom _{\hat{A}} (R^{0} (\Delta S) (\mathcal{F}),\O _{\hat{x}})\neq 0.
\]
Since $R^{i} (\Delta S) (\mathcal{F})=0$ for $i>0$ we have
\[
\Hom _{\hat{A}} (R^{0} (\Delta S) (\mathcal{F}),\O _{\hat{x}})=\Hom _{D
(\hat{A})} (\mathbf{R} (\Delta S) (\mathcal{F}),\O _{\hat{x}})
\]
where $\O _{\hat{x}}$ is regarded as a complex concentrated in degree
0. Now $\O _{\hat{x}}=\mathbf{R} (\Delta S) (L_{\hat{x}})$ since
$\mathbf{R}S (L_{\hat{x}})=\O _{\hat{x}}[2]$ and $\mathbf{R}\Delta (\O
_{\hat{x}}[2])=\O _{\hat{x}}$. So
\begin{align*}
\Hom _{D (\hat{A})} (\mathbf{R} (\Delta S) (\mathcal{F}),\O
_{\hat{x}})&=\Hom _{D (\hat{A})} (\mathbf{R} (\Delta S)
(\mathcal{F}),\mathbf{R} (\Delta S) (L_{\hat{x}}))\\
&=\Hom _{D ({A})} (L_{\hat{x}},\mathcal{F})\\
&=\Hom _{A} (L_{\hat{x}},\mathcal{F})
\end{align*}
using the fact that $\mathbf{R} (\Delta S)$ is an anti-equivalence of
categories. This gives us $\Hom (L_{\hat{x}},\mathcal{F})\neq 0$ which
implies $\Hom (L_{\hat{x}},\E )=H^{0} (\E \otimes L_{\hat{x}}^{\vee })\neq
0$. This in turn implies (by Proposition 4.18 of \cite{Mukai-Semi}) that
$H^{0} (\E ^{\vee }\otimes L_{\hat{x}})=\Hom (\E ,L_{\hat{x}})\neq 0$ from
which we get $\Hom (\mathcal{F},L_{\hat{x}})\neq 0$. Therefore we get a
(necessarily non-constant) endomorphism $\mathcal{F}\to L_{\hat{x}}\to
\mathcal{F}$ which contradicts the simplicity of $\mathcal{F}$.\qed

\subsection{Mukai stability}

We can regard the three equivalent conditions in Theorem~\ref{thm: Mukai's
thm} as giving a different stability condition for bundles. We call this
condition Mukai stability and we extend it to principal $G^{\cnums }$
bundles:

\begin{defn}\label{defn: Mukai stability}
We say a semi-stable holomorphic bundle $\E $ on $A$ is \emph{Mukai-stable}
if there exists a framing $\tau :\E \to \O _{p_{0}}$ such that $\Ker (\tau
)$ is simple, i.e. $\End (\mathcal{\Ker (\tau )})=\cnums $. We say that a
holomorphic $GL (n,\cnums )$, $SL (n,\cnums )$, or $Sp (n,\cnums )$ bundle
is Mukai stable if the associated vector bundle (induced by the standard
representation) is Mukai stable.
\end{defn}

Using the facts that $\Hilb ^{n} (\hat{A})$ is a fine moduli space and
$\mathbf{R} (\Delta {S})$ is an equivalence of categories, and applying
Theorem~\ref{thm: Mukai's thm}, we get the following (c.f. Theorem 2.20 of
\cite{Mukai-Senai}):

\begin{theorem}\label{thm: moduli of Mukai stable bundles is Hilb}
Let $\til{M} (n,0,0)$ be the space of Mukai-stable bundles $\E$ on $A$ with
$ch (\E )= (n,0,0)$, then $\til{M} (n,0,0)$ is a fine moduli space and the
functor $\mathbf{R} (\Delta S)$ applied to $\Ker (\tau )$ induces an
isomorphism $\til{M} (n,0,0)\cong \Hilb ^{n} (\hat{A})$. Moreover the map
$\til{M} (n,0,0)\to M (n,0,0)$ induced by sending $\E $ to its
$s$-equivalence class fits into the following commutative diagram
\[
\begin{diagram}
\til{M} (n,0,0)&\rTo&M (n,0,0)\\
\dTo&&\dTo\\
\Hilb ^{n} (\hat{A})&\rTo&\Sym ^{n} (\hat{A})
\end{diagram}
\]
where the vertical arrows are isomorphisms induced by $\mathbf{R}
(\Delta S)$ applied to $\Ker (\tau )$ and $\E$ respectively. In
particular, $\til{M} (n,0,0)$ is a resolution of singularities of $M
(n,0,0)$.
\end{theorem}

This theorem provides an answer to Question \ref{question: is the
resolution a moduli space?} in the case when $G$ is the (non-semi-simple)
group $U (n)$. That is, it shows that the moduli space $M_{U (n)} (A)=M
(n,0,0)$ of flat $U (n)$ connections on $A$ has a hyperk\"ahler resolution
given by the moduli space $\til{M} (n,0,0)$ of Mukai stable $U (n)^{\cnums
}$-bundles (\ie $Gl (n,\cnums )$-bundles).

We now use this theorem to analyze the Mukai-stable moduli spaces in the
$Sp (n)$ and $SU (n)$ cases.

\begin{defn}\label{defn: tilMsu and tilMsp}
Let $\til{M}_{SU (n)} (A)\subset \til{M} (n,0,0)$ be the subset of Mukai
stable bundles that arise as the associated vector bundles of principal
holomorphic $SL (n,\cnums )=SU (n)^{\cnums }$ bundles. For the case of $Sp
(n)$, let $\til{M}_{Sp (n)} (A)\subset \til{M} (2n,0,0)$ be the
\emph{closure} of the subset of Mukai stable bundles that arise as the
associated vector bundles of principal holomorphic $Sp (n,\cnums )=Sp
(n)^{\cnums }$ bundles.
\end{defn}

We now wish to show (as the notation suggests) that
$\til{M}_{SU (n)} (A)$ and $\til{M}_{Sp (n)} (A)$ are hyperk\"ahler
resolutions of $M_{SU (n)} (A)$ and $M_{Sp (n)} (A)$.

We first treat the $SU (n)$ case. A holomorphic vector bundle $\E $ is the
associated bundle of a holomorphic principal $SL (n,\cnums )$ bundle if
and only if $\det\E \cong \mathcal{O}$. There fore $\til{M}_{SU (n)}
(A)\subset \til{M} (n,0,0)$ is the subset of bundles with trivial
determinant. Since the determinant of a bundle is constant in its
$s$-equivalence class, $\til{M}_{SU (n)} (A)$ is the preimage of the subset
of bundles in $M (n,0,0)$ which have trivial determinant. This set is just
$M_{SU (n)} (A)$ which under the isomorphism $M (n,0,0)\cong \Sym ^{n}
(\hat{A})$ is the fiber over $p_{0}$ of the sum map $\Sym ^{n} (\hat{A})\to
\hat{A}$. Thus we recover the \emph{ad hoc} hyperk\"ahler resolution
constructed in Section \ref{sec: G-bundles on E and A}:

\begin{thm}\label{thm: tilMsu is the Mukai stable moduli space}
The moduli space of Mukai stable $SL (n,\cnums )$ bundles $\til{M}_{SU (n)}
(A)$ is smooth, holomorphic symplectic, and isomorphic to the generalized
Kummer variety $KA_{n-1}$. The natural map $\til{M}_{SU (n)} (A)\to M_{SU
(n)} (A)$ is a hyperk\"ahler resolution.
\end{thm}

We now turn to the $Sp (n)$ case. A holomorphic bundle $\E $ is the
associated bundle of a holomorphic principal $Sp (n,\cnums )$ bundle if and
only if there is an isomorphism
\[
\phi :\E \to \E ^{\vee }
\]
such that $\phi ^{\vee }=-\phi $ (\ie a symplectic form). We apply the
Fourier-Mukai transform to translate this into a condition for $\til{M}_{Sp
(n)} (A)\subset \til{M} (2n,0,0)\cong \Hilb ^{2n} (\hat{A})$. The result is
the following.

\begin{thm}\label{thm: mukai stable Sp bundles is hyperkahler and birat to HilbX}
The moduli space of Mukai stable $Sp (n,\cnums )$ bundles $\til{M}_{Sp
(n)} (A)$ (Definition \ref{defn: tilMsu and tilMsp}) is smooth, holomorphic
symplectic, and birationally equivalent to $ \operatorname{Hilb}^{n} (X)$
where $X$ is the Kummer $K3$ surface associated to $A$. The natural map
$\til{M}_{Sp (n)} (A)\to M_{Sp (n)} (A)$ is a hyperk\"ahler resolution.
\end{thm}

\begin{remark}\label{rem: Msp has two different resolutions (maybe)}
Note that strictly speaking, Definition \ref{defn: tilMsu and tilMsp} is
not consistent with our earlier \emph{ad hoc }definition of $\til{M}_{Sp
(n)} (A)$ as $\Hilb ^{n} (X)$ since Theorem \ref{thm: mukai stable Sp
bundles is hyperkahler and birat to HilbX} only asserts that they are
birationally equivalent. However, from the point of view of this paper, the
distinction is not very important---we have compared the Hodge numbers of
the resolution to the stringy Hodge numbers and so we can work with either
of the desingularizations since the Hodge numbers of birationally
equivalent hyperk\"ahler manifolds are the same. As we remarked earlier, it
is believed (but has not been proved) that birationally equivalent
hyperk\"ahler manifolds are actually deformation equivalent and hence
diffeomorphic, c.f. \cite{Huybrechts}.

We also note that the subset of Mukai stable bundles that admit a
symplectic form is not closed in $\til{M} (2n,0,0)$. Thus the points of
$\til{M}_{Sp (n)} (A)$ parameterize not only symplectic Mukai stable
bundles but also Mukai stable bundles with degenerate symplectic forms that
occur as the limits of symplectic Mukai stable bundles. We will give
examples of such bundles in the course of the proof of the theorem (see
Remark~\ref{rem: example non-symplectic bundle in tilMsp}).
\end{remark}

\textsc{Proof of Theorem \ref{thm: mukai stable Sp bundles is hyperkahler
and birat to HilbX}:} Since $\E $ is Mukai stable, by Theorem \ref{thm:
Mukai's thm}, $\mathbf{R}S (\E )=\hat{\E }[2]=\mathcal{O}_{\hat{Z}}[2]$
where $\hat{Z}\subset \hat{A}$ is a length $2n$, 0 dimensional subscheme of
$\hat{A}$. In the sequel, we drop the hats and just write $Z\subset A$. Let
$(-1): {A}\to  {A}$ be the involution $ {x}\mapsto - {x}$ and
note that $(-1)^{*}$ induces an involution on $\Hilb ^{2n} ( {A})$. We
use $- {Z}$ to denote the image of $ {Z} $ under $(-1)_{*}$ so that
we also have $(-1)^{*}\mathcal{O}_{ {Z}}=\mathcal{O}_{- {Z}}$. There
is a natural equivalence of functors \cite{Mukai-Duality}
\[
\mathbf{R} (S\Delta )= (-1)^{*}\mathbf{R} (\Delta S)[2]
\]
and so 
\begin{align*}
\mathbf{R}S (\E ^{\vee })&=\mathbf{R} (S\Delta ) (\E )\\
&= (-1)^{*}\mathbf{R}\Delta (\mathcal{O}_{ {Z}}[2])\\
&= (-1)^{*}\sheafext ^{2} (\mathcal{O}_{ {Z}},\mathcal{O})\\
&=\sheafext ^{2} (\mathcal{O}_{- {Z}},\mathcal{O}).
\end{align*}

For an arbitrary 0-dimensional subscheme $W\subset A$, the sheaf $\sheafext
^{2} (\mathcal{O}_{W},\mathcal{O})$ is not necessarily isomorphic to
$\mathcal{O}_{W}$ because it may fail to be the structure sheaf of a
subscheme. However, if $\sheafext ^{2} (\mathcal{O}_{W},\mathcal{O})$
\emph{is} the structure sheaf of a subscheme, then $\sheafext ^{2}
(\mathcal{O}_{W},\mathcal{O})$ is isomorphic to $\mathcal{O}_{W}$ (although
not canonically!); this assertion will be proved in the course of the
proof of Lemma \ref{lem: MtilSp is the fixed point set}. In the case at
hand, since $\E \cong \E ^{\vee }$ we have $\mathcal{O}_{Z}\cong \sheafext
^{2} (\mathcal{O}_{-Z},\mathcal{O})\cong \mathcal{O}_{-Z}$ and so $Z=
-Z$. In particular, $\til{M}_{Sp (n)} (A)$ is contained in the fixed point
set of the action of $(-1)^{*}$ on $\Hilb ^{2n} ( {A})$. More precisely, we
have the following Lemma:

\begin{lemma}\label{lem: MtilSp is the fixed point set}
$\til{M}_{Sp (n)} (A)$ is a connected component of the fixed locus of
$(-1)^{*}$ acting on $\Hilb ^{2n} (A)$.
\end{lemma}

We defer the proof of the lemma to the appendix. The theorem follows easily
from the lemma: First, since $\Hilb ^{2n} ( {A})$ is smooth, the components
of the fixed point set of the involution $(-1)^{*}$ are
smooth. Furthermore, since $(-1)^{*}$ preserves the holomorphic symplectic
form on $\Hilb ^{2n} ( {A})$, the fixed components are also holomorphic
symplectic. Finally, the fixed component that we claim is $\til{M}_{Sp (n)}
(A)$ lies over the subset of $\Sym ^{2n} ( {A})$ consisting of $n$ pairs of
points of the form $\{ {x},- {x} \}$. This set is naturally identified with
$\Sym ^{n} ( {A}/\pm 1)$. Thus $\til{M}_{Sp (n)} (A)$ is birational to
$\Sym ^{n} ( {A}/\pm 1)$ which is birational to $\Hilb ^{n} (X)$.\qed

\appendix \section{Miscellaneous details}\label{appendix}

Here we provide the details that were suppressed in the main discourse.

\subsection{The basic example for $D_{4}$.}

\begin{thm}\label{thm: Mspin(8) has a C8/-1 point}
Let $W$ and $\Lambda $ be the Weyl group and coroot lattice for $Spin
(8)$. There exists a point of $(A\otimes \Lambda )/W$ locally modeled on
$\cnums ^{8}/\pm 1$.
\end{thm}

\textsc{Proof:} We first derive a useful description of $A\otimes \Lambda
$. The coroot lattice $\Lambda $ is the sublattice of $\znums ^{4}$
generated by the simple coroots $e_{1}-e_{2}$, $e_{2}-e_{3}$,
$e_{3}-e_{4}$, and $e_{3}+e_{4}$. $W$ is generated by the reflections
through the planes perpendicular to the simple coroots.

Thus one can easily see that 
\[
\Lambda =\{\sum _{i=1}^{4}a_{i}e_{i}\in \znums ^{4}\text{ : }\sum
a_{i}\equiv 0 \bmod 2 \}
\]
and 
\[
W=S_{4}\ltimes \{\pm 1\}^{3}\subset S_{4}\ltimes \{\pm 1\}^{4}
\]
where the action of $W$ on $\Lambda $ is the restriction of the action on
$\znums ^{4}$ given by permuting the factors and multiplying by $-1$ on
some even number of factors. The elements $\sum a_{i}e_{i}$ with
$a_{i}\equiv 0 \bmod 2$ form a sublattice of $\Lambda $ giving us the exact
sequence:
\[
0\to (2\znums )^{4}\to \Lambda \to (\znums /2)^{3}\to 0
\]
where $(\znums /2)^{3}\subset (\znums /2)^{4}$ is the kernel of the sum map
$(\znums /2)^{4}\to \znums /2$. Noting that $A\otimes (2\znums )^{4}\cong
A\otimes \znums ^{4}= A^{4}$ we apply $A\otimes (\cdot )$ to the sequence
and examine the Tor sequence to arrive at:
\[
0\to \Tor _{1} (A,(\znums /2)^{3})\to A^{4}\to A\otimes \Lambda \to 0.
\]

The subgroup $\Tor _{1} (A,(\znums /2)^{3})\subset A^{4}$ is concretely
given as
\[
\Tor _{1} (A,(\znums /2)^{3})=\{(\tau _{1},\dots ,\tau _{4})\in A^{4}\text{
: $2\tau _{i}=0$, $\sum \tau _{i}=0$} \}.
\]

Thus we can regard elements of $A\otimes \Lambda $ as orbits of points in
$A^{4} $ by translation by the finite number of elements in $\Tor _{1}
(A,(\znums /2)^{3})$. The action of $W=S_{4}\ltimes \{\pm 1\}^{3}$ on
$A\otimes \Lambda $ is induced from the natural action on $A^{4}$. 

Now we choose 3 distinct, non-zero 2-torsion points of $A$ that sum to 0;
that is, let $\tau _{1},\tau _{2},\tau _{3}\in A$ be such that $2\tau
_{i}=0$, $\tau _{i}\neq \tau _{j}\neq 0$ for all $i\neq j$, and $\tau
_{1}+\tau _{2}+\tau _{3}=0$. Note that $(0,\tau _{1},\tau _{2},\tau
_{3})\in \Tor _{1} (A,(\znums /2)^{3})$.

We then choose ``square-roots'' of the $\tau _{i}$'s; that is,
elements $\tau _{i}/2\in A$, such that $2 (\tau _{i}/2)=\tau _{i}$.

Let 
\[
p= (0,\tau _{1}/2,\tau _{2}/2,\tau _{3}/2)\in A^{4}.
\]

In $A\otimes \Lambda $, $-p$ is equivalent to $p$ since
\begin{align*}
(0,-\tau _{1}/2,-\tau _{2}/2,-\tau _{3}/2)&=(0, 
\tau _{1}-\tau _{1}/2,\tau _{2}-\tau _{2}/2,\tau _{3}-\tau _{3}/2)\\
&= (0,\tau _{1}/2,\tau _{2}/2,\tau _{3}/2).
\end{align*}

Thus the subgroup $\{\pm 1 \}\subset S_{4}\ltimes \{\pm 1\}^{3}\subset
S_{4}\ltimes \{\pm 1\}^{4} $ with generator $Id\times (-1,-1,-1,-1)$ fixes
$p$ in $A\otimes \Lambda $. In fact, the stabilizer of $p $ in $A\otimes
\Lambda $ is exactly $\{\pm 1 \}$: The stabilizer does not contain
non-trivial permutations because $0$, $\tau _{1}/2$, $\tau _{2}/2$, and
$\tau _{3}/2$ are distinct. No other subgroup of $\{\pm 1\}^{3}$ stabilizes
$p$ since $\tau _{i}/2\neq -\tau _{i}/2$ and we must have at least two
$-1$'s acting.  Thus, the local model of the quotient of $A\otimes \Lambda
$ by $W$ near $p$ is $\cnums ^{8}/\pm 1$ (where $\pm 1$ acts non-trivially
on all factors) as asserted by Theorem~\ref{thm: Mspin(8) has a C8/-1
point}. \qed

\begin{remark}
If we replace $A$ with $E$ in the above discussions, we see that in
$M_{Spin (8)} (E)= (E\otimes \Lambda) /W$ there is a point modeled on
$\cnums ^{4}/\pm 1$. Looijenga's theorem (Theorem~\ref{thm: looijengas
thm}) tells us that $M_{Spin (8)} (E)$ is in fact $\cnums \P (1,1,1,1,2)$
which has a unique singular point (modeled on $\cnums ^{4}/\pm 1$). In an
elliptic curve $E$, there are exactly 3 non-zero 2-torsion points and so
the choice of the $\tau _{i}$ is unique (up to permutation) and the choice
of the square roots $\tau _{i}/2$ is unique up to addition by a 2-torsion
point. It is then easily checked that the $W$ orbit of the point $(0,\tau
_{1}/2,\tau _{2}/2,\tau _{3}/2)\in A\otimes \Lambda $ is unique and
corresponds to the predicted singular point in $\cnums \P (1,1,1,1,2)$. In
$A$, there are many choices for the $\tau _{i}$'s and so there are multiple
points in $M_{Spin (8)}^{0} (A)$ where a crepant resolution does not exist
locally.
\end{remark}

\subsection{The basic example for $B_{3}$.}

\begin{thm}\label{thm: Mspin(7) has a C8/-1 point}
Let $W$ and $\Lambda $ be the Weyl group and coroot lattice for $Spin
(7)$. There exists a point of $(A\otimes \Lambda )/W$ locally modeled on
$\cnums ^{6}/\pm 1$.
\end{thm}
\textsc{Proof:} The coroot lattice of $B_{n}$ in general is given as the
sublattice of $\znums ^{n}$ defined by the condition that the sum of the
coordinates should be even. That is, there is an exact sequence
\begin{equation}\label{eqn: exact seq LambdaBn-->Zn-->Z/2}
\begin{diagram}
0&\rTo& \Lambda _{B_{n}}&\rTo& \znums ^{n}&\rTo^{sum}&\znums /2&\rTo& 0.
\end{diagram}
\end{equation}

The Weyl group $W=W_{B_{n}}$ is the same as the Weyl group of the $C_{n}$
coroot system: $W_{C_{n}}\cong S_{n}\ltimes \{\pm 1 \}^{n}$. The action of
$W$ on $\Lambda _{B_{n}}$ is induced by the action on $\znums ^{n}\cong
\Lambda _{C_{n}}$ which is given by permuting the factors and multiplying
by $\pm 1$ on each coordinate. We remark that for the \emph{root} lattices,
the situation is exactly reversed, $\Lambda ^{*}_{C_{n}}\subset \Lambda
^{*}_{B_{n}}\cong \znums ^{n}$. These facts are easily seen by examining
the simple roots (see for example the tables in Appendix C of
\cite{Knapp}). The simple roots of the $B_{n}$ root system are given by
\[
\{e_{1}-e_{2},e_{2}-e_{3},\dots ,e_{n-1}-e_{n},e_{n} \}
\]
which span the full lattice $\znums ^{n}=\sum _{i}e_{i}\znums $, whereas the
simple roots of $C_{n}$ are given by
\[
\{e_{1}-e_{2},e_{2}-e_{3},\dots ,e_{n-1}-e_{n},2e_{n} \}
\]
which spans the kernel of the (mod 2) sum $\znums ^{n}\to \znums /2$. To
obtain the coroot system from the root system, one replaces each root
$\alpha $ with its coroot (see for example page 496 of
\cite{Fulton-Harris})
\[
\alpha '=\frac{2}{||\alpha ||^{2}}\alpha .
\]

It is then clear that the coroot system of $B_{n}$ is isomorphic to the
root system of $C_{n}$ and vice versa since
$(e_{i}-e_{i+1})'=e_{i}-e_{i+1}$, $e_{n}'=2e_{n}$ , and $(2e_{n})'=e_{n}$.

As in the $D_{4} $ case (indeed, $D_{n}$ in general) we can regard $\Lambda
_{B_{n}}$ as a quotient of $(2\znums )^{n}$:
\[
0\to (2\znums )^{n}\to \Lambda _{B_{n}}\to (\znums /2)^{n-1}\to 0
\]
where $(\znums /2)^{n-1}$ is the kernel of the sum map $(\znums /2)^{n}\to
\znums /2$. Thus, as in the $D_{4}$ case, we have a sequence
\[
0\to \Tor _{1} (A,(\znums /2)^{n-1})\to A^{n}\to A\otimes \Lambda
_{B_{n}}\to 0.
\]

Specializing to $n=3$, we can express the subgroup $\Tor _{1} (A,(\znums
/2)^{n-1})\subset A^{3}$ as
\[
\{(\tau _{1},\tau _{2},\tau _{3})\in A^{3}\text{ : $2\tau _{i}=0$ and $\sum
\tau _{i}=0$} \}.
\]
As in the $D_{4} $ case, choose $(\tau _{1},\tau _{2},\tau _{3})\in \Tor
_{1} (A,(\znums /2)^{n-1})$ such that $\tau _{i}\neq \tau _{j}\neq 0$ for
all $i\neq j$ and choose elements $\tau _{i}/2\in A$ such that $2 (\tau
_{i}/2)=\tau _{i}$. Define 
\[
p= (\tau _{1}/2,\tau _{2}/2,\tau _{3}/2)\in A^{3}/\Tor _{1} (A,(\znums
/2)^{n-1}).
\]
By essentially the same argument as in the $D_{4} $ case, the stabilizer of
$p$ is $\{\pm 1 \}$. Therefore the local model of the image of $p$ in
$(A\otimes \Lambda _{B_{3}})/W$ is $\cnums ^{6}/\{\pm 1 \}$.\qed

\subsection{The proof of Lemma \ref{lem: MtilSp is the fixed point set}}
Let $F^{0}\subset \Hilb ^{2n} (A)$ be the locus of subschemes consisting of
$2n $ distinct points of the form $\{p_{1},-p_{1},\dots ,p_{n},-p_{n}
\}$. Let $F$ be the closure of $F^{0}$. Note that $F$ has dimension $2n$ an
is a connected component of the fixed locus of $(-1)^{*}$ and is hence
smooth. Let $S=S_{n}\subset \Hilb ^{2n} (A)$ be the locus of subschemes
whose Fourier-Mukai transforms admit a symplectic form so by definition
$\til{M}_{Sp (n)} (A)=\overline{S}$. We will prove that $\overline{S}=F$.

Let $H_{k}\subset \Hilb ^{2k} (A)$ be the locus of subschemes supported at
the origin and let $H_{k}'\subset \Hilb ^{2k} (A)$ be the locus of
subschemes supported at the two-torsion points. By Theorem~1.13 of
\cite{Nakajima}, $\dim H_{k}\leq 2k-1$ for $k\geq  1$. It follows that
$\dim H'_{k}\leq 2k-1$. In fact, a component of $H_{k}'$
parameterizing subschemes supported at precisely $l$ of the two-torsion
points has dimension $2k-l$.

Now $S$ has a stratification 
\[
S=\bigcup _{k\geq 0}S^{0}_{n-k}\times (H'_{k}\cap S_{k}),
\]
where $S^{0}_{n-k}\subset S_{n-k}\subset \Hilb ^{2 (n-k)} (A)$
parameterizes subschemes whose support is \emph{disjoint} from the
two-torsion points. Since $\dim (S^{0}_{n-k})=2 (n-k)$, we see that all
strata of $S$ have dimension less than or equal to $2n-1$, except for
$S^{0}_{n}$ which is irreducible of dimension $2n$: In fact, there is a
quasi-finite map from an open subset of $\Hilb ^{n} (A)$ onto $S^{0}_{n}$.
Note that $F^{0}$ is a dense open subset of $S^{0}_{n}$. Indeed, if
$Z=\{p_{1},-p_{1},\dots ,p_{n},-p_{n} \}$ is a subscheme corresponding to a
point in $F^{0}$, then the Fourier-Mukai transform of $\O _{Z}$ is the
direct sum $L_{p_{1}}\oplus L_{-p_{1}}\oplus \dots \oplus L_{p_{n}}\oplus
L_{-p_{n}}$ of degree zero line bundles which has an obvious symplectic
form.

We claim that all the strata of $S$ are contained in the closure
$\overline{S^{0}_{n}}$ and hence $S=\overline{F^{0}}=F$ which proves the
lemma. This claim is equivalent to showing that $\dim _{z} (S)=2n$ at every
point $z\in S$. By factoring $S$ locally near $z$, this reduces to showing
that $\dim _{z} (S_{k})=2k$ at each $z\in S_{k}\cap H_{k}$. Since
$n$ is arbitrary, we prove this for $k=n$.

To prove this claim, we begin by examining the Fourier-Mukai transform of
the condition $-\phi =\phi^{\vee } $. We get $-\mathbf{R}S (\phi
)=\mathbf{R}S (\phi ^{\vee })=\mathbf{R} (S\Delta ) (\phi )=
(-1)^{*}\mathbf{R} (\Delta S) (\phi )[2]$ and so the following diagram
commutes:
\begin{equation}\label{diagram: fourier transform of symplectic condition}
\begin{diagram}
\mathcal{O}_{ {Z}}&\rTo^{-R^{2}S (\phi )}&\sheafext ^{2} (\mathcal{O}_{- {Z}},\mathcal{O})\\
\dTo<{(-1)^{*}}&&\dTo>{(-1)^{*}}\\
\mathcal{O}_{- {Z}}&\rTo^{R^{0} (\Delta S) (\phi ) }&\sheafext ^{2} (\mathcal{O}_{ {Z}},\mathcal{O})\\
\end{diagram}
\end{equation}

Let $V_{Z}=H^{0} (\mathcal{O}_{ {Z}})$ so by Serre duality $V_{Z}^{\vee
}=\operatorname{Ext}^{2} (\mathcal{O}_{ {Z}},\mathcal{O})$. Applying the
global sections functor $\Gamma $ to the above diagram and writing $\Phi $
for $\mathbf{R} (\Gamma S) (\phi )$ we get
\[
\begin{diagram}
V_{Z}&\rTo^{-\Phi }&V_{-Z}^{\vee }\\
\dTo<{(-1)^{*}}&&\dTo>{(-1)^{*}}\\
V_{-Z}&\rTo^{\Phi ^{\vee }}&V_{Z}^{\vee }
\end{diagram}
\]
so that in particular, $(-1)^{*}\circ \Phi $ is a symplectic form on $V_{Z}$.

Suppose now that $z\in S_{n}\cap H_{n}$. That is, $z$ corresponds to a
subscheme $Z\subset A$ with $\operatorname{Supp} ( {Z})=p_{0}$
such that we have the Diagram~(\ref{diagram: fourier transform of
symplectic condition}). The sheaf $\mathcal{O}_{ {Z}}$ is then determined
by the corresponding module over the local ring $
\hat{\mathcal{O}}_{p_{0}}\cong \cnums [[x,y]]$. This has a concrete
description in terms of linear algebra (c.f. Nakajima \cite{Nakajima}
section 1.2), namely $\mathcal{O}_{ {Z}}$ is determined by the actions of
$x$ and $y$ on $V_{Z}$. That is, $\mathcal{O}_{ {Z}}$ is determined by a
pair of nilpotent endomorphisms $M_{x},M_{y}\in \End (V_{Z})$ that
commute. Conversely, suppose $M_{x}$ and $M_{y}$ are any pair of commuting,
nilpotent, $2n\times 2n $ dimensional complex matrices. Then the action of
$M_{x}$ and $M_{y}$ give $\cnums ^{2n}$ the structure of a finite length
module over $\cnums [[x,y]]$; this module will be of the form $\cnums
[[x,y]]/\mathcal{I}_{Z}$ (and hence correspond to a point in $H_{0}$) if
and only if there exists a vector $v\in \cnums ^{2n}$ such that the vectors
$\{M^{i}_{x}M^{j}_{y} (v) \}_{i,j\geq 0}$ span $\cnums ^{2n}$. In this
case, $\mathcal{I}_{Z}=\{f\in \cnums [[x,y]]: f (M_{x},M_{y})=0 \}$ and the
matrices are unique up to simultaneous conjugation.

More generally, if $(M_{x},M_{y})$ are a pair of (not necessarily
nilpotent) commuting matrices satisfying the above spanning condition, then
they determine (uniquely up to simultaneous conjugation) an ideal
$\mathcal{I}\subset \cnums [x,y]$ of finite length and hence a 0
dimensional subscheme of $\cnums ^{2}$.

Note that $H_{n}\subset \Hilb ^{2n} (A)$ has a neighborhood $\nu (H_{n})$,
open in the analytic topology, parameterizing subschemes whose support is
contained in an open $\epsilon $-polydisc about $p_{0}$. $\nu (H_{n})$ is
isomorphic to the corresponding neighborhood of $H_{n}\subset \Hilb ^{2n}
(\cnums ^{2})$ whose points are given by (equivalence classes of) commuting
matrices $(M_{x},M_{y})$ whose eigenvalues have modulus less than $\epsilon
$. Under this identification, Diagram~(\ref{diagram: fourier transform of
symplectic condition}) for $A$ corresponds to the same diagram for $\cnums
^{2}$. We may therefore replace $A$ and the subschemes $S_{n}$ and $H_{n}$
of $\Hilb ^{2n} (A)$ by $\cnums ^{2}$ and the corresponding subschemes of
$\Hilb ^{2n} (\cnums ^{2})$.

If the sheaf $\mathcal{O}_{Z}$ is given by the pair $(M_{x},M_{y})$, then
the sheaf $\mathcal{O}_{- {Z}}$ is then clearly given by the pair
$(-M_{x},-M_{y})$. We can also determine the matrices corresponding to the
sheaf $\sheafext ^{2} (\mathcal{O}_{ {Z}},\mathcal{O})$. This sheaf must be
given by the pair $(M_{x}^{t},M_{y}^{t})$ where $M_{x}^{t},M_{y}^{t}$
denote the transpose matrices. This follows from the uniqueness of the
dualizing functor for modules over a local ring (\cite{HartRD} pg. 275):
both $\mathcal{O}_{Z}\mapsto \sheafext ^{2} (\mathcal{O}_{Z},\mathcal{O})$
and $(M_{x},M_{y})\mapsto (M_{x}^{t},M_{y}^{t})$ satisfy the conditions of
a dualizing functor and so they must coincide\footnote{One sees from this
description that the sheaf $\sheafext ^{2} (\mathcal{O}_{Z},\mathcal{O})$
is isomorphic (non-canonically) to $\mathcal{O}_{Z}$ as long as $\sheafext
^{2} (\mathcal{O}_{Z},\mathcal{O})$ is of the form $\mathcal{O/I}$. This
follows since the ideals $\{f\in \cnums [[x,y]]: f (M_{x},M_{y})=0 \}$ and
$\{f\in \cnums [[x,y]]: f (M_{x}^{t},M_{y}^{t})=0 \}$ coincide. However, it
can happen that $(M_{x}^{t},M_{y}^{t})$ fails the spanning condition even
if it is met by $(M_{x},M_{y})$. For example, $\sheafext ^{2} (\mathcal{O}/
(x^{2},xy,y^{2}),\mathcal{O})$ is not isomorphic to the structure sheaf of
a subscheme. In this example, the pair
$(M_{x},M_{y}) = \left( \left(\begin{smallmatrix}0&0&0\\
1&0&0\\
0&0&0
\end{smallmatrix} \right),\left(\begin{smallmatrix}0&0&0\\
0&0&0\\
1&0&0
\end{smallmatrix} \right) \right)$ satisfies the spanning condition, but the
transpose pair $(M^{t}_{x},M_{y}^{t})$ does not.}.  In this language, the
existence of Diagram~(\ref{diagram: fourier transform of symplectic
condition}) means that there is a skew-symmetric, invertible matrix $\Phi
=-\Phi ^{t}$ such that $\Phi M_{\bullet }\Phi ^{-1}=-M_{\bullet }^{t}$ for
$\bullet =x,y$.

Note that whenever $(M_{\bullet },\Phi )$ satisfy $\Phi M_{\bullet }\Phi
^{-1}=-M^{t}_{\bullet } $, then $(P^{-1}M_{\bullet }P,P^{t}\Phi P)$ satisfy
the same equation. Since $\Phi $ is skew-symmetric and invertible, there
exists a $P$ so that $P^{t}\Phi P=J$ where $J=\left(\begin{smallmatrix}0&I\\
-I&0\end{smallmatrix} \right)$ is the standard symplectic form. Thus we may
assume that $(M_{x},M_{y})$ satisfy $JM_{\bullet }+M^{t}_{\bullet }J=0$,
\ie $(M_{x},M_{y})\in \mathfrak{sp} (n)\oplus \mathfrak{sp} (n)$. This
equation is preserved under $M_{\bullet }\mapsto g^{-1}M_{\bullet }g$ if
and only if $g^{t}Jg=J$, \ie $g\in Sp (n,\cnums )$.

We conclude that $S\subset \Hilb ^{2n} (\cnums ^{2})$ is the image of an
open set $\mathcal{C}^{0}$ of the ``commuting variety''
\[
\mathcal{C}^{0}\subset \mathcal{C}=\{(g_{1},g_{2})\in \mathfrak{sp}\oplus
\mathfrak{sp}:[g_{1},g_{2}]=0\}
\]
under a morphism whose fibers are the simultaneous $Sp (n,\cnums )$
conjugacy classes in $\mathcal{C}^{0}$. The open set $\mathcal{C}^{0}$ is
the set of $(g_{1},g_{2})\in \mathcal{C}$ satisfying the spanning
condition; that is, there exists a vector $v\in \cnums ^{2n}$ such that
$\{g_{1}^{i}g_{2}^{j}v \}_{i,j\geq 0}$ spans $\cnums ^{2n}$. By Theorem~A
of Richardson \cite{Richardson}, $\mathcal{C}$ is irreducible. It follows
that $S\subset \Hilb ^{2n} (\cnums ^{2})$, the image of $\mathcal{C}^{0}$,
is also irreducible, hence everywhere of dimension at least $2n$, and so
the same must hold for the original $S\subset \Hilb ^{2n} (A)$.\qed
 
\begin{remark}\label{rem: example non-symplectic bundle in tilMsp}
We have shown that $\overline{S}=F$ but in general $S\neq F$---the
symplectic form can degenerate in a family of Mukai stable symplectic
bundles. For example, let $\E $ be the Fourier-Mukai transform of $\O _{Z}$
where $Z $ is a subscheme supported at $p_{0}\in A$ defined by the ideal
$\mathcal{I}= (x,y)^{3}\cup (y^{2}-xy,x^{2}-xy)$. A matrix representation
of this subscheme is given by the pair
\[
M_{x}=\left(\begin{matrix}0&1&0&0\\
0&0&1&0\\
0&0&0&0\\
0&0&0&0
\end{matrix} \right),
M_{y}=\left(\begin{matrix}0&1&0&0\\
0&0&1&0\\
0&0&0&0\\
0&0&1&0
\end{matrix} \right)
\]

One can check by hand that there does not exist an \emph{invertible} matrix
$\Phi $ such that $\Phi =-\Phi ^{t}$ and $\Phi M_{\bullet }=-M^{t}_{\bullet
}\Phi $ for $\bullet =x,y$. Thus $\E$ is not symplectic.

Another example with a slightly different flavor is as follows. Let $\E$ be
the Fourier-Mukai transform of $\O / (\mathcal{I}_{p_{0}})^{2}$. As we
showed earlier in the footnote, $\O / (\mathcal{I}_{p_{0}})^{2}$ is
invariant under $(-1)^{*}$ but it is not isomorphic to $\sheafext ^{2} (\O
/ (\mathcal{I}_{p_{0}})^{2},\O )$. This implies that $\E \not\cong \E
^{\vee }$. Consequently, any bundle of the form $(L\otimes \E)\oplus
(L^{-1}\otimes \E)$ is Mukai stable but cannot have a non-degenerate
symplectic form (here $L$ is a degree zero line bundle that is \emph{not}
two-torsion). Examples of this type occur in codimension 4; in fact, one
can prove that in general, the components of the locus of bundles in
$\til{M}_{Sp (n)} (A)$ with degenerate symplectic forms have codimension at
least 4.
\end{remark}


\end{document}